\documentclass[11pt]{article}

\usepackage{amstex,amscd}
\usepackage{theorem}

\title{Non-Hermitian Yang-Mills connections.}
\author{D. Kaledin, M. Verbitsky.}

\def\cdot{{\:\raisebox{3.5pt}{\text{\circle*{1.5}}}}}

\newcommand{\arrow}{\longrightarrow}
\newcommand{\Tw}{\operatorname{Tw}}

\newcommand{\id}{\operatorname{id}}
\newcommand{\cp}{{{\Bbb CP}^1}}
\newcommand{\6}{\partial}
\newcommand{\E}{{\cal B}}
\newcommand{\EE}{{\cal E}}
\newcommand{\bE}{\overline{\E}^*}
\newcommand{\bnabla}{\overline{\nabla}}
\newcommand{\K}{{\cal K}}
\newcommand{\Ker}{\operatorname{Ker}}
\newcommand{\pr}{\operatorname{pr}}
\newcommand{\tr}{\operatorname{tr}}
\newcommand{\Id}{\operatorname{Id}}
\newcommand{\C}{{\Bbb C}}
\newcommand{\R}{{\Bbb R}}
\newcommand{\Z}{{\Bbb Z}}
\newcommand{\A}{{\cal A}}
\newcommand{\F}{{\cal F}}
\newcommand{\W}{{\cal W}}
\newcommand{\M}{{\cal M}^s}
\newcommand{\X}{{\cal X}}
\newcommand{\bM}{\overline{\cal M}^s}
\newcommand{\gM}{{{\cal M}^{gs}}}
\newcommand{\G}{{\cal G}}
\newcommand{\CC}{{\cal C}}
\newcommand{\HH}{{\cal H}}
\newcommand{\ex}{{\operatorname{ex}}}
\newcommand{\conj}{\overline{\phantom{\E}}}
\newcommand{\wt}{\widetilde}
\newcommand{\End}{\operatorname{{\cal E}\!{\it nd}}}
\newcommand{\Hom}{\operatorname{{\cal H}\!{\it om}}}
\newcommand{\Maps}{\operatorname{Maps}}
\newcommand{\Aut}{\operatorname{Aut}}
\newcommand{\Sec}{\operatorname{Sec}}
\newcommand{\Var}{{\frak V}{\frak a}{\frak r}}
\newcommand{\Sets}{\operatorname{Sets}}
\newcommand{\topp}{{\operatorname{top}}}
\newcommand{\fib}{{\operatorname{fib}}}
\newcommand{\rank}{\operatorname{rank}}

\newcommand{\emp}{\bf}
\newcommand{\bpi}{\overline{\pi}}
\newcommand{\eps}{\varepsilon}
\newcommand{\ups}{\upsilon}

\newcommand{\su}{{SU(2)}}
\newcommand{\ssu}{\frak{su}(2)}
\newcommand{\inv}{{\operatorname{inv}}}
\newcommand{\h}{{\Bbb H}}
\newcommand{\calo}{{\cal O}}

\newcommand{\proof}{\par\noindent{\bf Proof.}\ }
\def\blacksquare{\hbox{\vrule width 4pt height 4pt depth 0pt}}
\def\endproof{\blacksquare

\medskip

}

\newtheorem{theorem}{Theorem}[section]

\newtheorem{lemma}[theorem]{Lemma}

\newtheorem{prop}[theorem]{Proposition}
\newtheorem{corr}[theorem]{Corollary}
\newtheorem{conjecture}[theorem]{Conjecture}

\theorembodyfont{\rm}

\newtheorem{rem}[theorem]{Remark}
\newtheorem{defn}[theorem]{Definition}
\newtheorem{question}[theorem]{Question}
\newtheorem{example}[theorem]{Example}

\makeatletter

\newcommand{\ps@verbit}{%
  \renewcommand{\@oddhead}{%
	  \scriptsize
	  {NHYM bundles}
	  \hfil\tiny {June 29, 1996}}
  \renewcommand{\@evenhead}{\@oddhead}
  \renewcommand{\@oddfoot}{\hfil\thepage\hfil}
  \renewcommand{\@evenfoot}{\@oddfoot}}
 
\newcommand{\restrict}[1]{{|_{{\phantom{|}\!\!}_{#1}}}}
 
\makeatother

\begin{document}

\begin{center}
{\bfseries \LARGE Non-Hermitian Yang-Mills connections.}\\
\vskip 1em%
{\large D. Kaledin, M. Verbitsky}\footnote{\tt 
kaledin$@$ium.ips.ras.ru, verbit$@$thelema.dnttm.rssi.ru} 
      \vskip 1.5em%
    {\large June 29, 1996} \\ 
\end{center}

\hfill

\hspace{0.1\linewidth}\begin{minipage}[t]{0.8\linewidth}
{\bf Abstract}.{
\small
We study Yang-Mills connections on holomorphic bundles
over complex K\"ahler manifolds of arbitrary dimension,
in the spirit of Hitchin's and Simpson's study of 
flat connections. The space of non-Hermitian Yang-Mills (NHYM)
connections has dimension twice the space of Hermitian Yang-Mills
connections, and is locally isomorphic to the complexification
of the space of Hermitian Yang-Mills connections (which is, by
Uhlenbeck and Yau, the same as the space of stable bundles).
Further, we study the NHYM connections over hyperk\"ahler manifolds.
We construct direct and inverse twistor transform from 
NHYM bundles on a hyperk\"ahler manifold to 
holomorphic bundles over its twistor space.
We study the stability and the modular properties of
holomorphic bundles over twistor spaces, and prove that work of Li and
Yau, giving the notion of stability for bundles over non-K\"ahler 
manifolds, can be applied to the twistors.  We identify
locally the following two spaces: the 
space of stable holomorphic bundles on a twistor
space of a hyperk\"ahler manifold
and the space of rational curves in the twistor space of the 
``Mukai dual'' hyperk\"ahler manifold. 

}
\end{minipage}

\hfill

\tableofcontents

\section{Introduction.}

\subsection{An overview}

In this paper we study non-Hermitian Yang-Mills (NHYM) connections 
on a complex vector bundle ${\cal B}$ over a K\"ahler manifold. 
By definition, a connection
$\nabla$ in ${\cal B}$  is Yang-Mills if its curvature $\Theta$
satisfies 
\begin{equation} \label{intro-Yang-Mills_Equation_}
\begin{cases}
\Lambda(\Theta)&=\text{const}\cdot\id\\
\Theta \phantom{\Lambda()}&\in \Lambda^{1,1}(M, \End({\cal B})), 
\end{cases}
\end{equation}
where $\Lambda$ is the standard Hodge operator,  
and $\Lambda^{1,1}(M, \End({\cal B})$ is
the space of $(1,1)$-forms with coefficients in 
$\End({\cal B})$ (see Definition \ref{_NHYM_Definition_} for details). 
This definition is standard \cite{UY}, \cite{Donaldson:surfa}.
However, usually $\nabla$ is assumed to be compatible
with some Hermitian metric in $\cal B$. This is why we use the term
``non-Hermitian Yang-Mills'' to denote Yang-Mills connections which are
not necessarily Hermitian.

An important analogy for our construction is the one with flat 
connections on a complex vector bundle ${\cal B}$. 
Recall that when $c_1({\cal B})
= c_2({\cal B}) = 0$, Hermitian Yang-Mills  
connections are flat (L\"ubcke's principle; see \cite{S}). 
The moduli of flat, but not necessary unitary bundles is 
a beautiful subject, well studied in literature (see, e. g. \cite{S2}).
This space has dimension twice the dimension of the moduli space of
unitary flat bundles and has a natural holomorphic symplectic form.
Also, generic part of the moduli of non-unitary flat bundles is
equipped with a holomorphic Lagrangian\footnote{Lagrangian with
respect to the holomorphic symplectic form.} 
fibration over the space of unitary flat connections.

When $c_1({\cal B})
= c_2({\cal B}) = 0$, the flat connections are in 
one-to-one correspondence with those holomorphic structures on ${\cal B}$
which make it a polystable\footnote{By {\em polystable} we will always
mean ``a direct sum of stable''. Throughout the paper, stability is
understood in the sense of Mumford-Takemoto.} holomorphic bundle
\cite{UY}, \cite{S}. For arbitrary bundle ${\cal B}$, a similar statement
holds if we replace ``flat unitary'' by ``Hermitian Yang-Mills''.
Thus, it is natural
to weaken the flatness assumption and consider instead all Hermitian
Yang-Mills connections.  The non-Hermitian Yang-Mills connections that
we define correspond then to connections that are flat but not
necessarily unitary. The basic properties listed above for
non-unitary flat bundles hold here as well. We show that
the moduli space of NHYM connections
has dimension twice the dimension of the moduli of
Hermitian Yang-Mills connections and is naturally equipped with 
a holomorphic symplectic form.
As in the case of flat bundles, generic part of 
the moduli of NHYM connections has 
a holomorphic Lagrangian fibration over the
space of Hermitian Yang-Mills connections.

Let us give a brief outline of the paper.  Fix a compact K\"ahler
manifold $M$ with a complex vector bundle ${\cal B}$. Let $\M$ be the set of
equivalence classes of NHYM connections on ${\cal B}$, and let $\M_0 \subset
\M$ be the subset of connections admitting a compatible Hermitian
metric. Both sets turn out to have natural structures of complex
analytic varieties. Recall that $\M_0$ is a moduli space of stable
holomorphic bundles (\cite{UY}; see also \ref{_Uhle-Yau_Theorem_}). 
After giving the relevant 
definitions, in Section 1 we study the structure of $\M$ in the
neighborhood of $\M_0$. We prove that $\dim\M = 2 \dim \M_0$ in a
neighborhood of $\M_0$. Moreover, we identify the ring of germs of
holomorphic functions on $\M$ near $\M_0$ with the ring of real-analytic
complex-valued functions on $\M_0$. Thus an open neighborhood $U
\supset \M$ is a complexification of $\M_0$ in the sense of Grauert. We
also construct a holomorphic $2$-form on $\M$ which is symplectic in a
neighborhood of $\M_0$. This picture is completely analogous to that
for the space of flat connections, studied by Hitchin, Simpson and
others (\cite{H},\cite{S2}).

For concrete examples and applications of our theory, we
consider the case of NHYM-connections over a hyperk\"ahler
manifold (Definition \ref{defn.hyperkahler}). 
In this case, it is natural to modify the NHYM
condition. Every hyperk\"ahler manifold is equipped with a quaternion
action in its tangent space. Since the group of unitary quaternions
is isomorphic to $SU(2)$, a hyperk\"ahler manifold has the group
$SU(2)$ acting on its tangent bundle.
Consider the corresponding action of $SU(2)$ on the
space $\Lambda^\cdot(M)$ of differential forms 
over a hyperk\"ahler manifold $M$.
Then all $SU(2)$-invariant $2$-forms satisfy
\eqref{intro-Yang-Mills_Equation_} (Lemma \ref{primitive}). 
Thus, if the
curvature $\Theta$ of a bundle $({\cal B}, \nabla)$ is 
$SU(2)$-invariant, ${\cal B}$ is NHYM. Converse is
{\it a priori} non-true when $\dim_\R M >4$: there are
$2$-forms satisfying \eqref{intro-Yang-Mills_Equation_}
which are not $SU(2)$-invariant. However, as the discussion
at the end of Section \ref{_autodu_Section_} shows, for bundles over
compact manifolds $SU(2)$-invariance of the curvature is
a good enough approximation of the NHYM property.

A connection in ${\cal B}$ is called {\bf autodual} 
if its curvature is $SU(2)$-invariant
(Definition \ref{_autodual_Definition_}).
For $\dim_\R M=4$, the autoduality, in the sense of our
definition, is equivalent to the anti-autoduality in the
sense of 4-dimensional Yang-Mills theory. 
Hermitian autodual bundles were studied at great length in
\cite{Vb}.\footnote{In \cite{Vb}, the term {\bf hyperholomorphic}
was used for ``Hermitian autodual''.}
Most of this paper is dedicated to the study of non-Hermitian
autodual bundles over compact hyperk\"ahler manifolds.

Consider the natural action of $SU(2)$ in the cohomology
of a compact hyperk\"ahler manifold (see the beginning of Section
\ref{_autodu_Section_}). Let ${\cal B}$ be a bundle with 
the first two Chern classes $c_1({\cal B})$, $c_2({\cal B})$
$SU(2)$-invariant. In \cite{Vb} we proved that every Hermitian
Yang-Mills connection in ${\cal B}$ is autodual. It is natural to 
conjecture that in such a bundle every NHYM connection is autodual.
In Theorem~\ref{_NHYM-are-autodu_Theorem_}, we prove a weaker form of this
statement: namely, in a neighbourhood of the space of Hermitian
Yang-Mills connections, every NHYM connection
is autodual, assuming the first two Chern classes
are $SU(2)$-invariant. This is done by constructing an explicit
parametrization of this neighbourhood (Proposition~\ref{series}).

Throughout the rest of this paper (starting from Section
\ref{_twistors_Section_}) we study algebro-geometrical
aspect of autodual connections. Two interdependent 
algebro-geometric interpretations of autoduality arise. 
Both of these itterpretations are related to the twistor
formalism, which harks back to the works of Penrose and
Salamon \cite{Sal}. Twistor contruction is explained in detail in
 Section \ref{_twistors_Section_}; here we give a brief outline
of this formalism.

Every hyperk\"ahler manifold $M$ has a whole 2-dimensional sphere
of integrable complex structures, called {\bf induced
complex structures}; these complex structures correspond bijectively
to $\R$-algebra embeddings from complex numbers to quaternions 
(see Definition \ref{_induced_co_str_Definition_}). We identify
this 2-dimensional sphere with $\C P^1$. Gluing all induced
complex structures together with the complex structure
in $\C P^1$, we obtain an almost complex structure on the
product $M \times \C P^1$ (Definition
\ref{twistor}). As proven by Salamon \cite{Sal},
this almost complex structure is integrable. The complex manifold 
obtained in this way is called {\bf the twistor space
for $M$}, denoted by $\Tw(M)$. Consider the natural projections 
\[ \sigma:\; \Tw(M)= M\times \C P^1 \arrow M, \ \ \ 
\pi:\; \Tw(M)= M\times \C P^1 \arrow \C P^1;\] the latter
map is holomorphic. The key statement to the twistor
transform is the following lemma.

\begin{lemma}\label{_autodua_(1,1)-on-twi-intro_Lemma_}
Let $({\cal B}, \nabla)$ be a bundle with a 
connection over a hyperk\"ahler
manifold $M$, and \[ (\sigma^*{\cal B}, \sigma^* \nabla)\]  
be the pullback
of $({\cal B}, \nabla)$ to the twistor space. Then $(\sigma^*{\cal B},
\sigma ^* \nabla)$ is holomorphic if and only if $({\cal B}, \nabla)$
is autodual.
\end{lemma}

\proof This is a restatement of Lemma \ref{_autodua_(1,1)-on-twi_Lemma_}.
\endproof

This gives a natural map from the space of autodual connections
on $M$ to the space of holomorphic bundles on $\Tw(M)$. We prove that this
map is injective, and describe its image explicitly.

For every points $x\in M$, the set $\sigma^{-1}(x)$ is a complex
analytic submanifold of the twistor space. The projection 
$\pi\restrict{\sigma^{-1}(x)}:\; \sigma^{-1}(x) \arrow \C P^1$
gives a canonical identification of $\sigma^{-1}(x)$ with 
$\C P^1$. The rational curve $\sigma^{-1}(x) \subset \Tw(M)$
is called {\bf a horisontal twistor line in $\Tw(M)$} (Definition
\ref{twistor}).

The following proposition provides an inverse
to the map given by Lemma
\ref{_autodua_(1,1)-on-twi-intro_Lemma_}.

\begin{prop}
Let $M$ be a hyperk\"ahler manifold, $\Tw(M)$ its twistor space and 
${\cal E}$ a holomorphic bundle over $\Tw(M)$. Then ${\cal E}$ 
comes as a pullback of an autodual bundle $({\cal B}, \nabla)$
if and only if restriction of ${\cal E}$ to all horisontal 
twistor lines is trivial as a holomorphic vector bundle.
Moreover, this autodual bundle is unique, up to equivalence.
\end{prop}

\proof This is a restatement of
Theorem \ref{_twisto_transfo_equiva_Theorem_}.
\endproof

We obtained an identification of the set of equivalence classes of
autodual bundles with a subset of the 
set of equivalence classes of bundles
over the twistor space. We would like to interpret this 
identification geometrically, as an identification of
certain moduli spaces. The autodual bundles are NHYM.
The set of equivalence classes of stable NHYM bundles is equipped
with a natural complex structure and is finite-dimensional,
as we prove in Section \ref{_NHYM_Section_}. This 
general construction is used to build
the moduli space of autodual bundles.
It remains to define the notion of stability
for holomorphic bundles over the twistor space and to construct the
corresponding moduli space. The usual (Mumford-Takemoto)
notion of stability does not work, because twistor spaces
are not K\"ahler.\footnote{Moreover, as can be easily shown,
the twistor space of a compact hyperk\"ahler manifold
admits no K\"ahler metric.}
We apply results of Li and Yau \cite{yl}, who define a notion
of stability for bundles over complex manifolds
equipped with a Hermitian metric satisfying 
a certain condition (see \eqref{yl}). The twistor space $\Tw(M)$ is 
isomorphic as a smooth manifold
to $\C P^1 \times M$ and as such is equipped with
the product metric. This metric is obviously Hermitian.
We check the condition of Li and Yau for twistor spaces
by computing the terms of \eqref{yl} explicitly. This
enables us to speak of stable and semistable bundles
over twistor spaces. 

Let ${\cal E}$ be a holomorphic bundle over $\Tw(M)$ obtained
as a pullback of an autodual bundle on $M$. We prove that
${\cal E}$ is semistable. This gives a holomorphic interpretation
of the moduli of autodual bundles on $M$. This is the first of our
algebro-geometric interpretations. The second interpretation involves
significantly more geometry, but yields a more explicit moduli
space. 

Let $M$ be a compact hyperk\"ahler manifold, and ${\cal B}$ a 
complex vector bundle with first two Chern classes invariant
under the natural $SU(2)$-action. Let $\widehat M$ be the moduli space
for the Hermitian Yang-Mills connections on ${\cal B}$.\footnote{Such
connections are always autodual, \cite{Vb}.}
Then $\widehat M$ is equipped with a natural hyperk\"ahler structure
(\cite{Vb}). The first result of this type was obtained by Mukai
\cite{_Mukai:K3_}
in the context of his duality between K3-surfaces; we use the
term ``Mukai dual'' for $\widehat M$ in this more general situation.

Let $X\stackrel{\pi}{\arrow} \C P^1$, 
$\widehat X\stackrel{\hat \pi}{\arrow} \C P^1$ be the twistor 
spaces for $M$, $\widehat M$, equipped with the natural holomorphic
projections to $\C P^1$. For an induced complex structure $L$
on the hyperk\"ahler manifold $M$, we denote by $(M, L)$
the space $M$ considered as a complex K\"ahler manifold with $L$
as a complex structure. Identifying the set of induced
complex structures with $\C P^1$, we consider $L$ as a point 
in $\C P^1$. Then, the complex manifold 
$(M, L)$ is canonically isomorphic to 
the pre-image $\pi^{-1}(L) \subset X$. 
By $i_L$ we denote the natural embedding
$(M, L) = \pi^{-1}(L) \stackrel{i_L}{\hookrightarrow} X$.
Let $B$ be a stable holomorphic bundle over $(M, L)$, with the
complex vector bundle ${\cal B}$ as underlying complex vector space.
In \cite{Vb}, we produce a canonical identification
between the moduli space of such stable bundles and the 
space $\widehat M$ of autodual connections in ${\cal B}$. Let 
${\cal F}_{B} = {i_L}_* {B}$ be the coherent sheaf 
direct image of ${B}$ under $i_L$. The moduli space of such
sheaves ${\cal F}_{B}$ is naturally identified with $\widehat X$
(Section \ref{_lines_Section_}; see also \cite{Vb}). 

Consider a holomorphic section $s$ of the map 
$\hat \pi:\; \widehat X \arrow \C P^1$, 
that is, a holomorphic embedding 
$s:\; \C P^1 \arrow \widehat X$ such that  $s\circ \hat \pi = \id$.
The image of such embedding is called 
{\bf a twistor line in $\widehat X$} (Section 
\ref{_twistors_Section_}).

Let ${\cal E}$ be a vector bundle over $X$ such that
the pullback ${i_L}^* {\cal E}$ is stable for all induced
complex structures $L\in \C P^1$. Such a bundle ${\cal E}$
is called {\bf fiberwise stable} (Definition \ref{fib.st}). 
{}From Lemma \ref{gen.st}
it follows that fiberwise stable bundles are also stable, in the
sense of Li--Yau. We restrict our attention to those bundles 
${\cal E}$ which are, as $C^\infty$-vector bundles, isomorphic
to $\sigma^*({\cal B})$, where ${\cal B}$ is our original
complex vector bundle on $M$.

Every fiberwise stable bundle ${\cal E}$ on $X$ gives a twistor
line $s_{\cal E}:\; \C P^1 \arrow \widehat X$ in $\widehat X$, where
$s_{\cal E}$ associates a sheaf ${i_L}_* {i_L}^* {\cal E}$
to $L\in \C P^1$. Since points of $\widehat X$ are identified with
isomorphism classes of such sheaves, the sheaf 
${i_L}_* {i_L}^* {\cal E}$ can be naturally considered as a point in
$\widehat X$.

Clearly, the 
moduli $St_f(X)$ space of fiberwise stable bundles is open in the
moduli $St(X)$ of stable bundles on $X$. 
This gives a complex structure on $St_f(X)$. The set $Sec(\widehat X)$
of twistor lines
in $\widehat X$ is equipped with a complex structure as a subset of the
Douady space of rational curves in $\widehat X$.\footnote{Douady spaces
are analogues of Chow schemes, defined in the complex-analytic
(as opposed to algebraic) setting.}
We constructed a holomorphic
map from $St_f(X)$ to $Sec(\widehat X)$. We prove that this
map is in fact an isomorphism of complex varieties
(Theorem \ref{iso}). The direct and inverse twistor transform
give a canonical identification between the 
moduli of autodual bundles
on $M$ and an open subset of the 
moduli of semi-stable bundles on $X$. Thus, we obtain an 
identification of the moduli of autodual bundles on $M$
and the space of twistor lines in $\widehat X$.

We must caution the reader that in this introduction we mostly
ignore the fact that all our constructions use different
notions of stability; thus, all identifications are valid
only locally in the subset where all the flavours of stability 
hold. The precise statements are given in Sections 
\ref{_twisto-tra_Section_}--\ref{_lines_Section_}.

\subsection{Contents} 

\begin{itemize}

\item
Here are the contents of our article. 

\item The Introduction is in two
parts: the first part explains the main ideas of this paper, and
the second gives an overview of its content, section by section.
These two parts of Introduction are independent.
The introduction is also formally independent from the main part and vice
versa. The reader who prefers rigorous discourse
might ignore the introduction and start reading
from Section \ref{_NHYM_Section_}.

\item Section \ref{_NHYM_Section_} contains the definition of 
NHYM (non-\-Her\-mi\-tian Yang-\-Mills) connection. We give the definition
of  
$(0,1)$-sta\-bi\-lity for NHYM connections and consider the
natural forgetful map 
\begin{equation} \label{_(0,1_stable_to_holo-Equation_}
   \pi:\; {\cal M}^s \arrow {\cal M}^s_0 
\end{equation}
from the space of $(0,1)$-stable NHYM-connections
to the moduli space of stable holomorphic bundles. The fiber of this map
is described explicitly through a power series and the 
Green operator (Proposition \ref{series}). This map is
also used to show that the moduli space of $(0,1)$-stable NHYM-connections
is correctly defined and finite-dimensional 
(Corollary \ref{_NHYM_finite-dim_Corollary_}).

Uhlenbeck--Yau theorem (Theorem \ref{_Uhle-Yau_Theorem_}; 
see also \cite{UY}) provides 
a compatible Hermitian Yang-Mils connection
for every stable holomorphic bundle. 
This gives a section 
${\cal M}^s_0 \stackrel i \hookrightarrow {\cal M}^s$
of the map \eqref{_(0,1_stable_to_holo-Equation_}. 
We study the structure of ${\cal M}^s$ in the neighbourhood
of $i\left({\cal M}^s_0\right)$, and prove that this neighbourhood
is isomorphic to the complexification of ${\cal M}^s_0$ in the
sense of Grauert (Proposition \ref{_complexi_Graue_Theorem_}).

\item In Section \ref{_autodu_Section_}, we recall the definition of
a hyperk\"ahler manifold and consider NHYM bundles over
a complex manifold with a hyperk\"ahler metric. We define autodual
bundles over hyperk\"ahler manifolds 
(Definition \ref{_autodual_Definition_})
and show that all autodual bundles are NHYM 
(Proposition \ref{_autodual_is_NHYM_Proposition_}). We cite the result of
\cite{Vb}, which shows that all Hermitian Yang-Mills connections
on a bundle ${\cal B}$ are autodual, if the first two Chern classes
of ${\cal B}$ satisfy a certain natural assumption ($SU(2)$-invariance;
see Theorem \ref{HYM.inv}).
We also prove that, for a NHYM connection $\nabla$ 
sufficiently close to Hermitian, $\nabla$ is autodual
(Theorem \ref{_NHYM-are-autodu_Theorem_}). 

\item Further on, we restrict our attention to autodual connection over
hy\-per\-k\"ah\-ler manifolds.

\item Section \ref{_twistors_Section_} gives a number of definition
and preliminary results from algebraic geometry of the twistor
spaces. We define the twistor space for an arbitrary hyperk\"ahler
manifold (Definition \ref{twistor}). The twistor space is a complex
manifold equipped with a holomorphic projection onto $\C P^1$. For most
hyperk\"ahler manifolds (including all compact ones), the
twistor space does not admit a K\"ahler metrics. This makes
it difficult to define stability for bundles over twistor
spaces. We overcome this difficulty by applying results
of Li and Yau (\cite{yl}). We consider the 
differential form which is an imaginary part of the natural
Hermitian metric on the twistor space. 
To apply \cite{yl}, we compute explicitly the
de Rham differential of this form
(Lemma \ref{_differe_of_Hermi_on_twistors_Lemma_}). 

\item In Section \ref{_twisto-tra_Section_} we define the direct and
inverse twistor transform relating autodual bundles over a
hyperk\"ahler manifold $M$ and holomorphic bundles over the corresponding
twistor space $X$. There is a map from the set of isomorphism
classes of autodual bundles on $M$ to the set of isomorphism classes
of holomorphic bundles on $X$ (Lemma \ref{_autodua_(1,1)-on-twi_Lemma_}).
We show that this map is an embedding and describe its
image explicitly (Theorem \ref{_twisto_transfo_equiva_Theorem_}).

\item In Section \ref{_stabi_of_twi_tra_Section_}, we 
consider holomorphic bundles on the twistor space
obtained as a result of a twistor transform. We
prove semistability of such bundles.
Thus, twistor transform is interpreted as a map between moduli
spaces.

\item In Section \ref{_lines_Section_} we return to the
study the algebro-geometric
properties of the twistor space. For a compact hyperk\"ahler manifold
$M$ and a stable holomorphic bundle $B$ on $M$, we consider the
space $\widehat M$ of deformations of $B$. When the first two Chern classes
of $B$ are $SU(2)$-invariant, the space $\widehat M$ has a natural
hyperk\"ahler structure; this space is called then {\bf Mukai dual}
to $M$. Let $X$, $\widehat X$ be the twistor spaces for $M$ and $\widehat M$.
We interpret the space of stable bundles on $X$ in terms of rational
curves on $\widehat X$ (Theorem \ref{iso}).

\item In Section \ref{_conje_Section_}, we relate a number of conjectures and open
questions from the geometry of NHYM and autodual bundles.  

\end{itemize}

\section{The general case.}
\label{_NHYM_Section_}
\subsection{Definition of NHYM connections}
Let $M$ be a K\"ahler manifold of dimension $n$ 
with the real valued K\"ahler form $\omega$. 
Consider a complex vector bundle $\E$
on $M$. Denote by $\A^n(\E)$ the bundle of smooth $\E$-valued $n$-forms 
on $\E$. Let 
$$
\A^n(\E) = \bigoplus_{i+j=n} \A^{i,j}(\E)
$$
be the Hodge type decomposition. The bundle $\End\E$ of endomorphisms of $\E$ 
is also a complex vector bundle. As usual, let 
$$
L:\A^{\cdot,\cdot}(\End\E) \to \A^{\cdot+1,\cdot+1}(\End\E)
$$
be the operator given by multiplication $\omega$. Let 
$$
\Lambda:\A^{\cdot+1,\cdot+1}(\End\E) \to \A^{\cdot,\cdot}(\End\E)
$$ 
be the adjoint operator with respect to the trace form on $\End\E$. 

\begin{defn}\label{_NHYM_Definition_}
A connection $\nabla:\E \to \A^1(\E)$ is called 
{\emp non-Hermitian Yang-Mills} (NHYM for short) if its curvature 
$R \in \A^2(\End\E)$ is of Hodge type $(1,1)$ and satisfies 
$$
\Lambda \circ R = c \Id
$$
for a certain constant $c \in \C$. 
\end{defn}

\begin{rem}
This terminology is perhaps unfortunate, in that a NHYM connection can 
(but need not) be Hermitian. We use the term for lack of better one. 
\end{rem}

To simplify exposition, we will always consider only NHYM connections 
with the constant $c=0$. 

Let $\nabla$ be a NHYM connection on $\E$. Since the $(0,2)$-component 
of its curvature vanishes, the $(0,1)$-component 
$\nabla^{0,1}:\E \to \A^{0,1}(\E)$ defines a holomorphic structure on 
$\E$. We will call this {\emp the holomorphic structure associated 
to $\nabla$}. 

Let $\bE$ be the dual to the complex-conjugate to the complex bundle
$\E$.  Every NHYM-connection $\nabla$ obviously induces a NHYM
connection $\nabla^*$ on the dual bundle $\E^*$. Let $\bnabla$ be the
connection on $\bE$ complex-conjugate to $\nabla^*$. The connection
$\bnabla$ is also obviously NHYM. We will call it {\emp the adjoint
connection to $\nabla$}. The holomorphic structure on $\bE$ associated to
$\bnabla$ will be called {\emp the adjoint holomorphic structure
associated to $\nabla$}. Note that the adjoint holomorphic structure
depends only on the $(1,0)$-part $\nabla^{1,0}$ of the connection
$\nabla$.

\subsection{Stability and moduli of NHYM connections.}
Fix a compact K\"ahler manifold $M$ and a complex vector bundle $\E$ on $M$. 
Consider the space $\A$ of all connections on $\E$ and let $\A_0$ be 
the subspace of NHYM connnections. The space $\A$ is a complex-analytic   
Banach manifold, and $\A_0 \subset \A$ is an analytic subspace of $\A$. 
Let $\G = \Maps(M,\Aut\E)$ be the complex Banach-Lie group of automorphisms 
of $\E$. The group $\G$ acts on $\A$ preserving the subset $\A_0$. 

In order to obtain a good moduli space for NHYM connections, we 
need to impose some stability conditions. 

\begin{defn} A NHYM connection $\nabla$ is called {\emp $(0,1)$-stable} if the 
bundle $\E$ with the associated holomorphic structure is a stable 
holomorphic bundle.
\end{defn}

Further on, we sometimes use the term {\it stable} to denote
$(0,1)$-stable connections.

\begin{rem}\label{too.strong}
This definition is sufficient for our present purposes. However, it is
unnaturally restrictive. See
\ref{_stabili_for_hyper_redu-Definition_} 
for a more natural definition. 
\end{rem}

Let $\A^s \subset \A_0$ be the open subset of $(0,1)$-stable NYHM connections and 
let 
$$
\M = \A_s / \G
$$
be the set of {\em equivalence classes} of $(0,1)$-stable NHYM connections on 
$\E$ endowed with the quotient topology. 

Choose a connection $\nabla \in \A_0, \nabla:\E \to \A^1(\E)$. Let 
$\nabla = \nabla^{1,0} + \nabla^{0,1}$ be the type decomposition and extend 
both components to differentials 
\begin{align*}
D:\A^{\cdot,0}(\End\E) &\to \A^{\cdot+1,0}(\End\E)\\
\overline{D}:\A^{0,\cdot} &\to \A^{0,\cdot+1}(\End\E)
\end{align*}
The tangent space $T_\nabla(\A)$ equals $T_\nabla(\A) =
\A^1(\End\E)$. The NHYM equations define a complex-analytic map
$$
YM:\A \to \A^{2,0}(\End\E) \oplus \A^{0,2} \oplus \A^0(\End\E).
$$
It is easy to see that the differential of $YM$ at the point $\nabla$ is 
given by 
$$
YM_\nabla = D + \overline{D} + \Lambda \nabla: \A^1(\End\E) \to 
\A^{2,0}(\End\E) \oplus \A^{0,2}(\End\E) \oplus \A^0(\End\E)
$$
On the other hand, the differential at $\nabla$ of the $\G$-action on 
$\A$ is given by 
$$
\nabla:\A^0(\End\E) \to \A^1(\End\E)
$$

\begin{defn}
The complex 
$$
0 \to \A^0(\End\E) \to \A^1(\End\E) \to 
\A^{2,0}(\End\E) \oplus \A^{0,2}(\End\E) \oplus \A^0(\End\E)
$$
is called {\emp the deformation complex} of the NHYM connection $\nabla$. 
\end{defn}

\begin{rem}
The deformation complex has a natural structure of a differential graded 
Lie algebra. 
\end{rem}

\begin{prop}
The deformation complex is elliptic.
\end{prop}

\proof
Indeed, the complex 
$$
0 \to \A^0(\End\E) \to \A^{0,1}(\End\E) \to \A^{0,2}(\End\E)
$$
is the Dolbeault complex for the holomorphic bundle $\End\E$ and is 
therefore elliptic. Hence it is enough to prove that 
\begin{equation}\label{kernel}
0 \to \A^{1,0}(\End\E) \to \A^{2,0}(\End\E) \oplus \A^0(\End\E)
\end{equation}
is elliptic. By Kodaira identity $\Lambda D = \sqrt{-1} D^*$ on
$\A^{1,0}(\End\E)$, and this complex is the same as
$$
0 \to D \oplus D^*:\A^{1,0}(\End\E) \to \A^{2,0}(\End\E) \oplus \A^0(\End\E),
$$
where $D^*$ is defined by means of the trace form on $\End\E$. This 
complex is obviously elliptic. 
\endproof

\begin{corr}
Let $\wt{G} \subset \G$ be the stabilizer of $\nabla \in \A$. Then 
\begin{enumerate}
\item $\wt{G}$ is a finite dimensional complex Lie group.
\item There exists a finite dimensional locally closed 
      complex-analytic Stein subspace $\wt{\M} \subset \A$ containing 
      $\nabla$ and invariant under $\wt{G}$ such that the natural projection 
      $$
      \wt{\M} / \wt{G} \to \M
      $$
      is an open embedding. 
\end{enumerate}
\end{corr}

\proof This is the standard application of the Luna's slice theorem, 
see \cite{Kod}.                    
\endproof 

\begin{corr} \label{_NHYM_finite-dim_Corollary_}
The topological space $\M$ has a natural structure of a comp\-lex\--analytic 
space. 
\end{corr}

\proof 
Indeed, since $\wt{\M}$ is Stein, the quotient $\wt{\M} / \wt{G}$ is a Stein 
complex-analytic space. Now by the standard argument (\cite{Kod}) the 
induced complex analytic charts on $\M$ glue together to give a 
complex-analytic structure on the whole of $\M$. 
\endproof 

\subsection[Hermitian Yang-Mills bundles and the theorem of 
Uhlenbeck--Yau]{Hermitian Yang-Mills bundles and the theorem of
\\ Uhlenbeck--Yau}

For every complex bundle $\E$ on $M$ denote by $\M_0(\E)$ the moduli
space of stable holomorphic structures on $\E$. Fix $\E$ and consider
the space $\M$ of $(0,1)$-stable NHYM connections on $\E$. Taking the
associated holomorphic structure defines a map $\pi:\M \to \M_0(\E)$.

\begin{lemma}
The map $\pi$ is holomorphic.
\end{lemma}

\proof
Clear.
\endproof

Since every complex vector bundle admits an Hermitian metric, the complex
vector bundles $\E$ and $\bE$ are isomorphic. Therefore the moduli spaces
$\M_0(\E)$ and $\M_0(\bE)$ are naturally identified. Denote the space
$\M_0(\E) = \M_0(\bE)$ simply by $\M_0$, and let $\bM_0$ be the
complex-conjugate space.

Consider the open subset $\gM \subset \M$ of $(0,1)$-stable NHYM connections
on $\E$ such that the adjoint connection on $\bE$ is also $(0,1)$-stable.
Taking the adjoint holomorphic structure defines a map $\bpi:\gM \to
\bM_0$. This map is also obviously holomorphic.

\begin{lemma} \label{_compa_Hermi_pi=barpi_Lemma_}
A NHYM connection $\nabla \in \gM$ satisfies $\pi(\nabla) = \bpi(\nabla)$ 
if and only if it admits a compatible Hermitian metric. 
\end{lemma}

\proof Indeed, $\pi(\nabla) = \bpi(\nabla)$ if and only if there exists 
an isomorphism $h:\E \to \bE$ sending $\nabla$ to $\bnabla$. This isomorphism 
defines an Hermitian metric on $\E$ compatible with the connection $\nabla$. 
\endproof 

To proceed further we need to recall the following fundamental theorem. 

\begin{theorem}[Uhlenbeck,Yau] \label{_Uhle-Yau_Theorem_}
Every stable holomorphic bundle $\E$ on a K\"ahler manifold $M$ admits
a unique Hermitian Yang-Mills connection $\nabla$.  Vice versa, every
holomorphic bundle admitting such a connection is polystable.
\end{theorem}

We will call such a metric {\emp a Uhlenbeck-Yau metric} for the
holomorphic bundle $\E$. 

\begin{corr}
Let $\M_u \in \gM$ be subset of equivalence classes of Hermitian
connections. The product map $\pi \times \bpi:\gM \to \M_0 \times
\bM_0$ identifies $\M_u$ with the diagonal in $\M_0 \times \bM_0$.
\end{corr}

\proof 
Clear. 
\endproof 

\begin{rem}
Note that the subset $\M_u \subset \gM$ is not complex-analytic, but 
only real-analytic. 
\end{rem}

\subsection{Moduli of NHYM connections as a complexification of
the moduli of stable bundles}

Let $\nabla \in \M_u \subset \M$ be an Hermitian Yang-Mills
connection, and let $\F_\nabla = \pi^{-1}(\pi(\nabla)) \subset \M$ be
the fiber of $\pi$ over $\nabla$. In order to study the map $\pi
\times \bpi:\gM \to \M_0 \times \bM_0$ in a neighborhood of $\nabla$, we
first study the restriction of the map $\bpi$ to $\F_\nabla$. We begin
with the following.

\begin{theorem}\label{zhopa}
Let $\overline{D}:\E \to \A^{0,1}(\E)$ be a representative in the
equivalence class $\pi(\nabla) \in \M_0$ of holomorphic structures on
$\E$ and let $D:\E \to \A^{1,0}(\E)$ be the operator adjoint to
$\overline{D}$ with respect to the Uhlenbeck-Yau metric.  The fiber
$\F_\nabla$ is isomorphic to the set of all $\End\E$-valued
$(1,0)$-forms $\theta$ satisfying
\begin{equation}\label{formula}
\begin{cases}
D\theta + \theta \wedge \theta &= 0,\\
D^*\theta &= 0.
\end{cases}
\end{equation}
\end{theorem}

\proof 
To define the desired isomorphism, choose for any equivalence class
$\nabla_1 \in \pi^{-1}(\pi(\nabla))$ a representative $\nabla_1 =
\nabla^{1,0}_1 + \overline{D}:\E \to \A^1(\E)$.  Every two
representatives must differ by a gauge transformation $g:\E \to \E$.
The map $g$ must preserve the holomorphic structure
$\overline{D}$. However, this holomorphic structure is by assumption
stable. Therefore $g = c\Id$ for $c \in \C$, and the operator
$\nabla^{1,0}_1:\E \to \A^{1,0}(\E)$ is defined uniquely by its class
in $\F_\nabla$. Take $\theta = \nabla^{1,0}_1 - D$; the
equations~\eqref{formula} follow directly from the definition of NHYM
connections.
\endproof

In order to apply this Theorem, note that by the second of the
equations~\eqref{formula} every NHYM connection $\nabla_1 \in
\F_\nabla$ defines a $D^*$-closed $\End\E$-valued $(1,0)$-form
$\theta$. Complex conjugation with respect to the Uhlenbeck-Yau metric
$h_\nabla$ identifies the space of $D^*$-closed $\End\E$-valued
$(1,0)$-forms with the space of $\overline{D}^*$-closed $(0,1)$-forms,
and it also identifies the respective cohomology spaces. But the
cohomology spaces of the Dolbeault complex $\A^{0,\cdot}(\End\E)$ with
respect to $\overline{D}$ and $\overline{D}^*$ are both equal to the
space of harmonic forms, hence naturally isomorphic.  Collecting all
this together, we define a map $\rho:\F_\nabla \to
\overline{H^1(M,\End\E)}$ by the rule
\begin{equation} \label{_from_NHYM_to_classes_Equation_}
\nabla_1 \mapsto \langle \text{ class of } \theta \text{ in } H^1(M,\End\E)
\rangle.
\end{equation}

\begin{prop}\label{kur}
The map $\rho$ is a closed embedding in a neighborhood of $\nabla \in
\M_u \in \F_\nabla$.
\end{prop}

This Proposition can be deduced directly from Theorem~\ref{zhopa}.
However, we prefer to prove a stronger statement. To formulate it,
consider the adjoint holomorphic structure $\bpi(\nabla)$ on the
complex bundle $\bE \cong \E$. Recall the following standard fact from
the deformation theory of holomorphic bundles.

\begin{theorem}
Let $\bar\6:\E \to \A^{0,1}(\E)$ be a stable holomorphic structure on
a complex Hermitian bundle $\E$.  There exists a neighborhood $U
\subset \M_0$ of $\bar\6$ such that every $\bar\6_1 \in U$ can be
represented uniquely by an operator $\bar\6 + \theta:\E \to
\A^{0,1}(\E)$ satisfying
\begin{equation}\label{kuranishi}
\begin{cases}
\bar\6\theta = \theta \wedge \theta\\
\bar\6^*\theta = 0
\end{cases}
\end{equation}
The {\emp Kuranishi map} $U \to H^1(M,\End\E)$ defined by 
$$
\bar\6_1 \to \langle \text{ class of } \theta \text{ in }
H^1_{\bar\6^*}(\A^{0,\cdot}(\End\E)) \cong H^1(M,\End\E) \rangle 
$$
is a locally closed embedding. 
\end{theorem}

\begin{corr}\label{compl}
The map $\bpi:\F_\nabla \to \bM_0$ is biholomophic in a neighborhood $V$ 
of $\nabla \in \M_u \subset \gM$, and the map $\rho:V \to H^1(M,\End\E)$ 
is the composition of $\bpi:\F_\nabla \to \bM_0$ and the Kuranishi map.
\end{corr}

\proof 
Complex conjugation sends equations~\eqref{formula} precisely to
\eqref{kuranishi}, and thus establishes a bijection between
neighborhoods of $\nabla \in \F_\nabla$ and $\bpi(\nabla) \in \bM_0$.
\endproof 

This statement, in turn, implies the following. 

\begin{prop} \label{_complexi_Graue_Theorem_}
The product map $\pi \times \bpi:\gM \to \M_0 \to \bM_0$ is biholomorphic on 
an open neighborhood $U$ of the subset $\M_u \in \gM$. 
\end{prop}

\proof
Consider both $\gM$ and $\M_0 \times \bM_0$ as spaces over $\M_0$. 
The map $\pi \times \bpi$ is a map over $\M_0$, and it is locally 
biholomorphic on every fiber of the natural projections $\gM \to \M_0$, 
$\M_0 \times \bM_0 \to \M_0$. 
\endproof 

Thus $\dim\M = 2\dim\M_0$, and an open neighborhood $U$ of the
subspace $\M_u \subset \gM$ is the complexification of $\M_0$ is the
sense of Grauert.

\subsection{Holomorphic symplectic form on the moduli of
NHYM bundles}

In order to construct a holomorphic symplectic $2$-form on $\M$, we
need to restrict our attention to a smooth open subset of $\M$.

\begin{defn}
A NHYM connection $\nabla$ is called {\emp smooth} if both $\nabla \in \M$ and 
$\pi(\nabla) \in \M_0$ are smooth points. 
\end{defn}

Let $\nabla \in \M$ be a smooth NHYM connection and denote by
$\CC_\nabla$ its deformation complex.  By construction the holomorphic
tangent space $T_\nabla(\M)$ is identified with a subspace of the
first cohomology space $H^1(\CC_\nabla)$. By definition $\pi(\nabla)$
is a smooth point, and the tangent space $T_{\pi(\nabla)}(\M_0)$ is a
subspace of the first cohomology space $H^1(M,\End\E)$ of $M$ with
coefficients in $\E$ equipped with the induced holomorphic structure.

Consider the natural projection map $\CC_\nabla \to
\A^{0,i}(M,\End\E)$ from $\CC_\nabla$ to the Dolbeault complex of the
bundle $\End\E$.  Denote by
$$
\pr:H^1(\CC_\nabla) \to H^1(M,\End\E)
$$ 
the induced map on the cohomology spaces and let 
$$
W = \Ker\pr \subset H^1(\CC_\nabla). 
$$

\begin{prop}\label{diagramma}
Let $d\pi:T_\nabla(\M) \to T_{\pi(\nabla)}(\M_0)$ be the differential
of the map $\pi:\M \to \M_0$ in the smooth point $\nabla \in
\M_0$. The diagram
$$
\begin{CD}
0 @>>> T_\nabla(\F_\nabla) @>>> T_\nabla(\M) @>d{\pi}>>   
				   T_{\pi(\nabla)} @>>> 0 \\
@VVV @VVV  @VVV   @VVV                        @VVV              \\
9 @>>> W                   @>>> H^1(\CC_\nabla) @>{\pr}>> 
				   H^1(M,\End\E)   @>>> 0
\end{CD}
$$
is commutative. 
\end{prop}

\proof 
Clear.
\endproof 

We first construct a symplectic form on the space $H^1(\CC_\nabla)$. 
To do this, we first identify the space $W \subset H^1(\CC_\nabla)$. 

\begin{lemma}
The space $W$ is naturally isomorphic to 
$H^{n-1}(M, \End\E \otimes \K)$, where $K$ 
is the canonical line bundle on $M$. 
\end{lemma}

\proof Indeed, the space $W$ is isomorphic to 
the space of $\End\E$-valued $(1,0)$-forms satisfying 
\begin{equation}\label{pipupa}
\begin{cases}
D\theta = 0\\
\Lambda \overline{D} \theta = 0
\end{cases}
\end{equation}
Consider the map 
$$
\bullet \wedge \omega^{n-1}:\A^{1,0}(\End\E) \to \A^{n,n-1}(\End\E).
$$
By Kodaira identities a form $\theta \in \A^1(M,\End\E)$ satisfies
\eqref{pipupa} if and only if $\theta \wedge \omega^{n-1}$ is
harmonic. Hence $\bullet \wedge \omega^{n-1}$ identifies $W$
with $H^{n-1}(M,\End\E \otimes \K)$.  
\endproof

Consider now two $\End\E$-valued $1$-forms $\theta_0,\theta^1 \in
\A^1(M,\End\E)$ and let
$$
\Omega(\theta_0,\theta_1) = \int_M 
\tr(\theta_0 \wedge \theta_1 \wedge \omega^{n-1}), 
$$
where $\tr$ is the trace map, $n = \dim M$ and $\omega$ is the K\"ahler 
form on $M$. 

\begin{lemma}
If $\theta_0 = \nabla g$ for some section $g \in \A^0(M,\End\E)$, then 
$$
\Omega(\theta_0,\theta_1) = 0
$$
for any $\theta_1 \in \A^1(M,\End\E)$ satisfying 
$\Lambda \nabla(\theta_1) = 0$.
\end{lemma}

\proof Indeed, 
\begin{multline*}
\Omega(\theta_0,\theta_1) = \int_M \tr( \nabla(\omega^{n-1}g)) 
\wedge \theta_1) = \\
= \int_M \tr(\omega^{n-1} \wedge \nabla(g\theta_1)) - 
\int_M \tr( g \omega^{n-1} \wedge \nabla(\theta_1)) = 
\int_M \tr(g \Lambda \nabla(\theta_1)) \omega^n = 0
\end{multline*}
\endproof

\begin{corr}\label{bububuj}
\begin{enumerate}
\item The form $\Omega$ defines a 
      complex $2$-form on the space $H^1(\CC_\nabla)$. 
\item The subspace $W$ is isotropic, and the induced pairing 
      $$
      (W \cong H^{n-1}(M,\End\E\otimes\K))\otimes H^1(M,\End\E) \to \C
      $$
      is non-degenerate. 
\end{enumerate}
\end{corr}

\proof 
The first statement is clear. It is easy to see that the pairing
induced by $\Omega$ is exactly the one defined by the Serre duality,
which proves the second statement.
\endproof  

Restricting to the subspace $T_\nabla(\M)$, we get a 
$2$-form $\Omega$ on $\M$. This form is obviously holomorphic. 

\begin{prop}\label{symplectic}
Assume that either $\nabla \in \M$ is Hermitian, or $T_\nabla(\M) =
H^1(\CC_\nabla)$.  Then the form $\Omega$ on $T_\nabla(\M)$ is
non-degenerate.  The map $\M \to \M_0$ is a Lagrangian fibration in
the neghrborhood of $\nabla$.
\end{prop}

\proof It is easy to see that $T_\nabla(\F_\nabla) \subset T_\nabla(\M)$ is
isotropic.  If the connection $\nabla$ is {\em unobstructed}, that is, the
embedding $T_\nabla(\M) \hookrightarrow H^1(\CC_\nabla)$ is actually an
isomorphism, then the statement follows from Corollary~\ref{bububuj}.

Suppose now that the inclusion $T_\nabla(\M) \subset H^1(\CC_\nabla)$
is proper.  By assumption the connection $\nabla$ is Hermitian in this
case, therefore the complex conjugation map $\conj:\A^{0,1}(M,\End\E)
\to \A^{1,0}(\End\E)$ is defined.  It is easy to see that this map
identifies $\Ker\pr$ with $\overline{H^1(M,\End\E)}$.  By
Corollary~\ref{compl} it also identifies $T_\nabla(\F_\nabla)$ with
$\overline{T_{\pi(\nabla)}}$.

The form $\Omega(\bullet,\overline{\bullet})$ is a non-degenerate
Hermitian form on $\Ker\pr$. Therefore its restriction to
$T_\nabla(M_\nabla)$ is also non-degenerate. 

The last statement now follows directly from Proposition~\ref{diagramma} and 
Corollary~\ref{bububuj}. 
\endproof 

\subsection{Local parametrization of the moduli of NHYM connections}

In the last part of this section we give a more explicit description
of the embedding $\F_\nabla \to H^1(\End\E)$ for an Hermitian
Yang-Mills connection $\nabla$ in the spirit of \cite{Vb}. This
description is of independent interest, and we will also use it in the
next section in the study of NHYM connections on hyperk\"ahler
manifolds.

Fix an Uhlenbeck-Yau metric on $\E$ compatible with $\nabla$ and let
$\Delta = DD^* + D^*D$ be the associated Laplace operator on
$\A^\cdot(\End\E)$. Let $G$ be the Green operator provided by the
Hodge theory. Recall that we have the Hodge decomposition
$$
\A^\cdot(\End\E) = \HH^\cdot(\End\E) \oplus \A^\cdot_\ex(\End\E)
$$
into the space of harmonic form $\HH^\cdot$ and its orthogonal complement 
$\A^\cdot_\ex$. This complement is further decomposed as 
$$
\A^\cdot_\ex(\End\E) = D\A^{\cdot-1}(\End\E) \oplus D^*\A^{\cdot+1}(\End\E)
$$
The composition $DD^*G$ is by definition the 
projection onto $D\A^{\cdot-1}(\End\E) \subset \A^\cdot(\End\E)$. 

Take now a small neighborhood $U \subset \F_\nabla$ of the Hermitian
connection $\nabla \in \F_\nabla$ and a NHYM connection $\nabla_1 \in
U$. Let $\theta = \nabla_1 - \nabla$, $\theta \in \A^{1,0}(M,\End\E)$
and let $K(\theta) \in \overline{H^{1,0}(M,\End\E)}$ be the associated
cohomology class. Shrinking $U$ if neccesary, we see that by
Proposition~\ref{kur} the connection $\nabla_1$ is uniquely determined
by the class $K(\theta)$.

Let $\theta_0 \in \A^{1,0}(M,\End\E)$ be the harmonic form
representing the class $K(\theta)$. Define by induction 
$$
\theta_n = D^*G \sum_{0 \leq k < n} \theta_k \wedge \theta_{n-1-k}.
$$

\begin{prop}\label{series}
Let $\nabla$ be an Hermitian Yang-Mills connection. 
There exists a neighborhood $V \subset \F_\nabla$ of 
$\nabla \in \F_\nabla$ such that 
for every $\nabla_1 \in V \subset \F_\nabla$ the series 
\begin{equation}\label{rjad}
\sum_{0 \leq k} \theta_k
\end{equation}
 converges 
to the form $\theta = \nabla_1 - \nabla$. 
\end{prop}

\proof 
The metric on $\End\E$ defines a norm $\|\bullet\|$ on
$\A^{\cdot,0}(\End\E)$. We can assume that
$\left\| \theta_0 \right\| < \eps$
for any fixed $\eps > 0$. Since the Hodge decomposition is orthogonal, 
\begin{multline*}
\left\|D \theta_n\right\| = \left\|DD^*G\left(\sum_{0 \leq k < n}
\theta_k \wedge \theta_{n-k}\right)\right\| \\ 
\leq \left\|\sum_{0 \leq k < n} \theta_k \wedge \theta_{n-k}\right\| 
\leq \sum_{0 \leq k < n}
\left\|\theta_k\right\| \cdot \left\|\theta_{n-k}\right\|.
\end{multline*}
Since $D:D^*(\A^{2,0}(\End\E)) \to \A^{2,0}(\End\E)$ is injective
and elliptic, there exists a constant $C > 0$ such that 
$$
\|Df\| > C\|f\|
$$ 
for all $f \in D^*(\A^{2,0}(\End\E))$. 
Let $a_n = \frac{(2n)!}{(n!)^2}$ be the Catalan numbers. By induction
$$
\left\|\theta_n\right\| < a_n \left| \frac{\eps}{C} \right|.
$$
Since $A(z) = \sum a_n z^n$ satisfies $A(z) = 1 + z(A(z))^2$, it equals 
$$
A(z) = \frac{1 - \sqrt{1-4z}}{2}
$$
and converges for $z < \frac{1}{4}$. Therefore the series~\eqref{rjad}
converges for $4\eps < C$. 

To prove that it converges to $\theta$, let $\chi_0 = \theta$ and let 
$$
\chi_n = \theta - \sum_{0 \leq k < n} \theta_k
$$
for $n \geq 1$.  Since both $\nabla$ and $\nabla + \theta$ are NHYM,
we have $D\chi_1 = \chi_0 \wedge \chi_0$. Therefore $\chi_0
\wedge \chi_0 \in \A^{2,0}(\End\E)$ and
$$
D \chi_1 = D D^* G (\chi_0 \wedge \chi_0) = \chi_0 \wedge \chi_0.
$$
By induction 
$$
D \chi_n = \chi_0 \wedge \chi_{n-1} + 
\sum_{0 \leq k \leq n} \chi_k \wedge \theta_{n - 1 - k},  
$$
and $\chi_n \in D^*(\A^{2,0}(M,\End\E))$ for all $n > 0$.  
Again by induction, 
$$
\left\|\chi_n\right\| < a_{n+1} \left(\frac{\eps}{C}\right)^n.
$$
Therefore $\chi_n \to 0$ if $4\eps < C$, which proves the Proposition. 
\endproof 

\begin{rem}
This Proposition can be strengthened somewhat. Namely, for any
harmonic $\End_\E$-valued $(1,0)$-form $\theta_0$ in a small
neighborhood of $0$ the series~\eqref{rjad} converges to a form
$\theta$. As follows from \cite{Vb}, the connection $\nabla + \theta$
is NHYM provided the following holds.
\begin{description}
\item[*] All the forms $\theta_p \wedge \theta_q \in \A^{2,0}(M,\End\E)$ 
	 lie in $\A^{2,0}_\ex(M,\End\E)$. 
\end{description}
This condition is also knows as vanishing of all of the so-called
Massey products $[\theta_0 \wedge \ldots \wedge \theta_0]$.
\end{rem}

\section{Autodual and NHYM connections in the hy\-per\-k\"ah\-ler case.}
\label{_autodu_Section_}

We now turn to the study of NHYM connections on hyperk\"ahler manifolds. 
First we recall the definitions and some general facts. 

\begin{defn}[\cite{Cal}]
\label{defn.hyperkahler}
A {\emp hyperk\"ahler manifold} is a Riemannian manifold $M$ equipped with 
two integrable almost complex structures $I,J$ which are parallel with 
respect to the Levi-Civita connection and satisfy 
$$
I \circ J = - J \circ I.
$$
\end{defn}

Let $M$ be a hyperk\"ahler manifold. The operators $I,J$ define an
action of the quaternion algebra $\h$ on the tangent bundle $TM$. This
action is also parallel. Every imaginary quaternion $a \in \h$
satisfying $a^2 = -1$ defines an almost complex structure on $M$.
This almost complex structure is parallel, hence integrable and
K\"ahler. 

\begin{defn} \label{_induced_co_str_Definition_}
A complex structure on $M$ corresponding to an imaginary quaternion $a \in \h$ 
with $a^2 = -1$ is said to be {\emp induced by $a$}. 
\end{defn}

For every such $a \in \h$ we will denote by $\omega_a$ the K\"ahler
form in the complex structure induced by $a$.  We will always assume
fixed a preferred complex K\"ahler structure $I$ on $M$.

Recall that every hyperk\"ahler manifold is equipped with a canonical 
holomorphic symplectic $2$-form $\Omega$. If $J,K \in \h$ satisfy $J^2
= -1$, $IJ=K$ then this form equals
$$
\Omega = \omega_J + \sqrt{-1}\omega_K.
$$

The group $U(\h)$ of all unitary quaternions is isomorphic to
$\su$. Thus every hyperk\"ahler manifold comes equipped with an action
of $\su$ on its tangent bundle, and, {\it a posteriori}, with an action of
its Lie algebra $\ssu$.  Extend these actions to the bundles
$\Lambda^\cdot$ of differential forms and let $\Lambda^\cdot_\inv
\subset \Lambda^\cdot$ be the subbundle of $\su$-invariant forms. The
$\su$-action does not commute with the de Rham differential. However,
it does commute with the Laplacian (see \cite{Vl}).  Therefore it
preserves the subspace of harmonic forms.  Identifying harmonic forms
with coholomogy classes, we get an action of $\su$ on the cohomology
spaces $H^\cdot(M)$.

Let $\Lambda:\Lambda^{\cdot+2} \to \Lambda^\cdot$ be the Hodge operator
associated to the K\"ahler metric on $M$. 

\begin{defn}
A differential form $\theta$ on $M$ is called {\emp primitive} if
$\Lambda\theta = 0$.
\end{defn}

\begin{lemma}\label{primitive}
\begin{enumerate}
\item All $\su$-invariant forms are primitive. All $\su$-in\-va\-ri\-ant
$2$-forms are of Hodge type $(1,1)$ for every one of the induced complex
structures on $M$. Vice versa, if a form is of type $(1,1)$ for all the
induced complex structures, it is $\su$-invariant.
\item The same statements hold for de Rham cohomology classes instead of
forms. 
\end{enumerate}
\end{lemma}

\proof See \cite{Vb}, Lemma 2.1 \endproof 

\begin{rem}
The converse is true for $\dim_\C M = 2$, but in higher dimensiona
there are primitive forms that are not $\su$-invariant.
\end{rem}

Consider a complex bundle $\E$ on $M$ and let $\nabla:\E \to \A^1(\E)$
be a connection on $\E$.

\begin{defn}\label{_autodual_Definition_}
The connection $\nabla$ is called {\emp autodual} if its curvature $R
\in \A^2(M,\End\E)$ is $\su$-inavariant.
\end{defn}

\begin{rem}
The terminology comes from the $4$-dimensional topology: autodual
connections on hyperk\"ahler surfaces are anti-selfdual in the usual
topological sense.
\end{rem}

We will call an autodual connection $\nabla$ {\emp $(0,1)$-stable} if its $(1,0)$-part 
defines a stable holomorphic structure on the bundle $\E$. Denote by 
$\M_\inv$ the set of equivalence classes of $(0,1)$-stable autodual connections on 
the bundle $\E$. 

\begin{rem}
Like in Remark~\ref{too.strong}, this $(0,1)$-stability condition may be too
restrictive.
\end{rem}

Let $\nabla$ be an autodual connection on $\E$. By Lemma~\ref{primitive}
for every $J \in \cp$ the $(0,1)$-component of the connection $\nabla$ with
respect to the complex structure induced by $J$ defines a holomorphic
structure on $\E$. We will call it {\emp the holomorphic structure induced
by $J$}.

\begin{prop} \label{_autodual_is_NHYM_Proposition_}
Let $M$ be a hyperk\"ahler manifold and let $\E$ be a complex bundle
on $M$. Every autodual connection $\nabla$ on $\E$ is NHYM.
\end{prop}

\proof 
Immediately follows from Lemma~\ref{primitive}. 
\endproof 

Therefore there exists a natural embedding $\M_\inv \hookrightarrow \M$
from $\M_\inv$ to the moduli space $\M$ of NHYM connections on $\E$. In the
rest of this section we give a partial description of the image of this
embedding.

We will use the following. 

\begin{theorem}\label{HYM.inv}
Assume that the first two Chern classes 
$$
c_1(\E),c_2(\E) \in H^\cdot(M)
$$ 
of the bundle $\E$ are $\su$-invariant. Then every Hermitian Yang-Mills
connection $\nabla$ on $\E$ is autodual.
\end{theorem}

\proof 
See \cite{Vb}. 
\endproof 

Therefore, if the Chern classes $c_1(\E),c_2(\E)$ are
$\su$-invariant, then the closed subset $\M_u \in \M$ of
Hermitian Yang-Mills connections lies in $\M_\inv$.

\begin{theorem} \label{_NHYM-are-autodu_Theorem_}
Let $M$ be a hyperk\"ahler manifold and let $\E$ be a complex vector
bundle on $M$ such that the Chern classes $c_1(\E),c_2(\E)$ are
$\su$-invariant.  Then subset $\M_\inv \subset \M$ of autodual
connections contains an open neighborhood of the subset $\M_u \subset \M$
of connections admitting a compatible Hermitian metric. 
\end{theorem}

\proof 
Let $\nabla$ be an Hermitian Yang-Mills connection. It is autodual by
Theorem~\ref{HYM.inv}, and it is enough to show that every $\nabla_1
\in \F_\nabla$ sufficiently close to $\nabla$ is also autodual.

Let $\theta = \nabla_1 - \nabla$ and let $\theta_n, n \geq 0$ be as in
\eqref{rjad}. By Propostion~\ref{series} we can assume that
$$
\theta = \sum_k \theta_k.
$$
It is enough to prove that $\overline{D}\theta$ is $\su$-invariant. We
will prove that $\overline{D}\theta_k$ is $\su$-invariant for all $k
\geq 0$. To do this, we use results of \cite{Vb}. 

Consider the operator 
$$
L^\Omega:\A^{\cdot,\cdot}(\End\E)\to\A^{\cdot+2,\cdot}(\End\E)
$$ 
given by multiplication by the canonical holomorphic $2$-form
$\Omega$ on $M$.  By \cite{Vl} the operator
$$
[L_\Omega,\Lambda]:\A^{\cdot,\cdot}(\End\E) \to \A^{\cdot+1,\cdot-1}(\End\E)
$$
coincides with the action of a nilpotent element in the Lie algebra $\ssu$.
Therefore it is a derivation with respect to the algebra structure on the
complex $\A^{\cdot,\cdot}(\End\E)$.  Moreover, a $(1,1)$-form $\alpha$ is
$\su$-invariant if and only if $[L_\Omega,\Lambda]\alpha = 0$.

Let $\6^J:\A^{\cdot,\cdot}(\End\E) \to \A^{\cdot+1,\cdot}(\End\E)$ 
be the commutator 
$$
\6^J = \left[\overline{D} , \left[L_\Omega,\Lambda\right]\right].
$$
This map again is a derivation with respect to the algebra structure
on $\A^{\cdot,\cdot}(\End\E)$, and it is enough to prove that
$\6^J\theta_n=0$ for all $n \geq 0$. Moreover, 
$$
\6^J = \left[\overline{D}, \left[L_\Omega,\Lambda\right]\right] =
\left[L_\Omega, \left[\overline{D}, \Lambda\right]\right] + \left[\Lambda,
\left[L_\Omega,\overline{D}\right]\right],
$$
The second term is zero since $\Omega$ is holomorphic, the first term is
$[L_\Omega,\sqrt{-1}D^*]$ by Kodaira identity.  Hence $D^*$ and $\6^J$
anticommute. Finally, the Laplacian $\Delta_J = \6_J\6_J^* + \6_J^*\6_J$ is
proportional by \cite{Vb} to the Laplacian $\Delta = DD^* +
D^*D$. Therefore the Laplacian $\Delta$ and the Green operator $G$ also
commute with $\6^J$.

Now, the form $\theta_0$ is by definition $\Delta$-harmonic.
Therefore it is also $\Delta_J$-harmonic, and $\6^J\theta_0 = 0$. To
prove that $\6^J\theta_n = 0$, use induction on $n$. By definition 
\begin{multline*}
\6^J\theta_n = \6^JD^*G\left(\sum_{0 \leq k < n} \theta_k \wedge
\theta_{n-1-k}\right) = \\
= -D^*G\left(\sum_{0 \leq k <n} \6^J\theta_k \wedge \theta_{n-1-k} + 
\theta_k \wedge \6^J\theta_{n-1-k} \right). 
\end{multline*}
The right hand side is zero by the inductive assumption. 

Thus all $\6^J\theta_k$ are zero, and all the $\overline{D}\theta_k$ are 
$\su$-invariant, which proves the Theorem. 
\endproof

\section{Stable bundles over twistor spaces.}
\label{_twistors_Section_}

\subsection{Introduction}
To further study autodual connections on a bundle $\E$ over a hyperk\"ahler
manifold $M$, we need to introduce the so-called ``twistor space'' $X$ for
$M$. This is a certain non-K\"ahler complex manifold associated to
$M$. Autodual connections give rise to holomorphic bundles on $X$ by means
of a construction known as ``twistor transform''. This construction turns
out to be essentially invertible, thus providing additional information on
the moduli space $\M_\inv$.
 
We develop the twistor transform machinery in the next section. In this
section we give the necessary preliminaries: the definition and some
properties of the twistor space $X$, and a discussion of the notion of
stability for holomorphic bundles over $X$. 

\subsection{Twistor spaces}
Let $M$ be a hyperk\"ahler manifold. Consider the product manifold $X = M
\times S^2$. Embed the sphere $S^2 \subset \h$ into the quaternion algebra
$\h$ as the subset of all quaternions $J$ with $J^2 = -1$. For every point
$x = m \times J \in X = M \times S^2$ the tangent space $T_xX$ is
canonically decomposed $T_xX = T_mM \oplus T_JS^2$. Identify $S^2 = \cp$
and let $I_J:T_JS^2 \to T_JS^2$ be the complex structure operator. Let
$I_m:T_mM \to T_mM$ be the complex structure on $M$ induced by $J \in S^2
\subset \h$.

The operator $I_x = I_m \oplus I_J:T_xX \to T_xX$ satisfies $I_x \circ I_x =
-1$. It depends smoothly on the point $x$, hence defines an almost complex
structure on $X$. This almost complex structure is known to be integrable
(see \cite{Sal}). 

\begin{defn}\label{twistor}
The complex manifold $\langle X, I_x \rangle$ is called {\emp the twistor
space} for the hyperk\"ahler manifold $M$. 
\end{defn}

By definition the twistor space comes equipped with projections $\sigma:X
\to M$, $\pi: X \to \cp$. The second projection is holomorphic.  For any
point $m \in M$ the section $\wt{m}:\cp \to X$ with image $m \times \cp
\subset X$ is also holomorphic. We will call this section $\wt{m}$ {\emp
the horizontal twistor line} corresponding to $m \in \M$. 

Let $\iota:\cp \to \cp$ be the real structure on $\cp$given by the
antipodal involution. Then the product map 
$$
\iota = \id \times \iota:X \to X
$$
defines a real structure on the complex manfiold $X$. The following
fundamental property of twistor spaces is proved, e.g., in \cite{HKLR}.

\begin{theorem}\label{inv}
Let $M$ be a hyperk\"ahler manifold and let $X$ be its twistor space. Then
a holomorphic section $\cp \to X$ of the natural projection $\pi:X \to \cp$
is a horizontal twistor line if and only if it commutes with natural real
structure $\iota:X \to X$.
\end{theorem}

Let $\Sec$ be the Douady moduli space of holomorphic sections $\cp \to X$
of the projection $\pi:X \to M$. Then conjugation by $\iota$ defines a real
structure on the complex-analytic space $\Sec$. Theorem~\ref{inv}
identifies the susbet of real points of $\Sec$ with the hyperk\"ahler
manifold $M$. We will call arbitrary holomorphic sections 
$\cp \to X$ {\emp twistor lines} in $X$.

\subsection{Li--Yau theorem}
\label{ss.twistor}
The twistor space, in general, does not admit a K\"ahler metric. In order
to obtain a good moduli space for holomorphic bundles on $X$ we use a
genralization of the notion of stability introduced by Li and Yau in
\cite{yl}. We reproduce here some of their results for the convenience of
the reader.

Let $X$ be an $n$-dimensional Riemannian complex manifold and let
$\sqrt{-1}\omega$ be the imaginary part of the metric on $X$. Thus $\omega$
is a real $(1,1)$-form.

Assume that the form $\omega$ satisfies the following condition. 
\begin{equation}\label{yl}
\omega^{n-2} \wedge d\omega = 0.
\end{equation}
For a closed real $2$-form $\eta$ let 
$$
\deg\eta = \int_X \omega^{n-1} \wedge \eta.
$$
The condition \eqref{yl} ensures that $\deg\eta$ depends only on the
cohomology class of $\eta$.  Thus it defines a degree functional
$\deg:H^2(X,\R) \to \R$. This functional allows one to repeat verbatim the
Mumford-Takemoto definitions of stable and semistable bundles in this
more general situation. Moreover, the Hermitian Yang-Mills equations also
carry over word-by-word. 

Yau and Li proved the following.

\begin{theorem}[\cite{yl}]
Let $X$ be a complex Riemannian manifold satisfying \eqref{yl}. Then every
stable holomorpic bundle $\E$ on $X$ admits a unique Hermitian Yang-Mills
connection $\nabla$. Vice versa, every bundle $\E$ admitting an Hermitian
Yang-Mills connection is polystable. 
\end{theorem}

Just like in the K\"ahler case, this Theorem allows one to construct a good
moduli space for holomorphic bundles on $X$. (See \cite{Kod}.)

\subsection{Li--Yau condition for twistor space}

The twistor space $X = M \times \cp$ is equipped with a natural Riemannian
metric, namely, the product of the metrics on $M$ and on $\cp$.  To apply
the Li-Yau theory to $X$, we need to check the condition \eqref{yl}. First,
we identify the form $\omega$.

Let $\omega = \omega_M + \omega_\cp$ be the decomposition associated with
the product decomposition $X = M \times \cp$. By the definition of the
complex structure on $X$, the form $\omega_\cp$ is the pullback
$\pi^*\omega$ of the usual K\"ahler form on $\cp$, while $\omega_M$ is a
certain linear combination of pullbacks of K\"ahler forms on $M$ associated
to different induced complex structures. 

Let $W \in \h$ be the $3$-dimensional subspace of imaginary
quaternions. For every $a \in W$ the metric $\langle \cdot,\cdot\rangle$ on
the hyperk\"ahler manifold $M$ defines a real closed $2$-form $\omega_a$ on
$M$ by the rule 
$$
\omega_a(\cdot,\cdot) = \langle \cdot, a \cdot \rangle.
$$
This construction is linear in $a$, hence defines an embedding $W
\hookrightarrow \A^2(M,\R)$. Let $\W$ be the trivial bundle on $\cp$ with
the fiber $W$. The embedding $W \hookrightarrow \A^2(M,\R)$ extends to an
embedding
$$
\W \hookrightarrow \pi_*\sigma^*\A^2(M,\R) \subset \pi_*\A^2(X).
$$
Since $X = M \times \cp$, the de Rham differential $d = d_X:\A^\cdot(X) \to
\A^{\cdot+1}(X)$ decomposes into the sum $d = d_M + d_\cp$. The
differential $d_\cp$ defines a flat connection on the bundle
$\pi_*\sigma^*\A^2(M,\R)$. The subbundle $\W \subset
\pi_*\sigma^*\A^2(M,\R)$ is flat with respect to $d_\cp$.

The space $W$ is equipped with an euclidian metric, thus $W = W^*$.  Since
we have an embedding $\cp = S^2 \hookrightarrow W$, the bundle $\W =
TW|_\cp = T^*W|_\cp$ decomposes orthogonally
$$
\W = \R \oplus \calo(-2)
$$
into the sum of the conormal and the cotangent bundles to $\cp \subset W$.
The conormal bundle is the trivial $1$-dimesional real bundle $\R$, and the
cotangent bundle is isomorphic to the complex vector bundle $\calo(-2)$ on
$\cp$. The connection $d_x|_\W$ induces the trivial connection on $\R$ and the
usual metric connection on $\calo(2)$. The embedding $\W \to
\pi_*\sigma\A^2(M,\R)$ decomposes then into a real $2$-form
$$
\omega \in \sigma^*\A^2(M,\R)
$$
and a complex $\calo(2)$-valued $2$-form 
$$
\Omega \in \sigma^*\A^2(M,\R) \otimes \pi^*\calo(2).
$$
The form $\Omega$ is holomorphic, while the form $\sqrt{-1}\omega_M$ is
precisely the imaginary part of the Hermitian metric on $X$.

Let now $\ups \in \A^{2,1}(X)$ be the $(2,1)$-form corresponding to the
holomorphic form
$$
\Omega \in \A^{2,0}(X,\pi^*\calo(2))
$$ 
under the identification $\calo(2) \cong \A^{0,1}(\cp)$ provided by the
metric on $\cp$.

\begin{lemma}\label{_differe_of_Hermi_on_twistors_Lemma_}
$$
d \omega = \ups + \overline{\ups}.
$$
\end{lemma}

\proof Indeed, $d\omega_\cp = 0$ and $d_M\omega_M = 0$, therefore it is
enough to compute $d_\cp\omega_M$. The bundle $\W \in \A^2(X)$ is invariant
under $d_\cp$, and $d$ induces the trivial connection $\nabla$ on
$\W$. Let
$$
d = \begin{bmatrix} \nabla_\R &\theta_{10}\\\theta_{01} &
\nabla_{\calo(2)} \end{bmatrix}
$$
be the decomposition of $\nabla$ with respect to $\W = \R \oplus
\calo(2)$. The connection $\nabla_\R$ is trivial, therefore
$\nabla_\R\omega = 0$. An easy computation shows that
$$
\theta_{10} \in \A^{1}(\cp,\calo(2))
$$
induces the isomorphism $\A^1(\cp) \cong \calo(2)$. Thus 
$$
d\omega_M = d_\cp\omega_M = \theta_{10} \wedge \Omega = \ups +
\overline{\ups}.
$$
\endproof 

Now we can prove that the twistor space $X$ satisfies the condition
\eqref{yl}. 

\begin{prop}
The canonical $(1,1)$-form $\omega$ on the twistor space $X$ satisfies 
$$
\omega^{n-1} \wedge d\omega = 0
$$
for $n = \dim M =\dim X - 1$.
\end{prop}

\proof 
Let $x = m \times J \in X$ be a point in $X$. Choose a local coordinate $z$
on $\cp$ near the point $J \in \cp$. In a neighborhood of $x \in X$ we have 
$$
\ups = f(z)\Omega \wedge d\bar z
$$
for some holomorphic function $f(z)$. Therefore 
$$
\omega^{n-1} \wedge d\omega = f(z)\omega^{n-1} \wedge \Omega \wedge
d\bar z + \overline{f(z)}\omega^{n-1} \wedge \overline{\Omega} \wedge dz.
$$
But $\omega^{n-1} \wedge \Omega$ and $\omega \wedge \overline{\Omega}$ are
both forms on the $n$-dimensional manifold $M$, of Hodge types $(n+1,n-1)$
and $(n-1,n+1)$. Hence both are zero. To prove the Proposition, it remains
to see that $\omega^{n-1} - \omega_M^{n-1}$ is divisible by $\omega_\cp$,
and $\omega_\cp \wedge d z = \omega_\cp \wedge d \bar z = 0$.
\endproof

\section{Twistor transform.}
\label{_twisto-tra_Section_}

\subsection{Twistor transform}

We now introduce the direct and inverse twistor transforms which relate
autodual bundles on the hyperk\"ahler manifold $M$ and holomorphic bundles
on its twistor space $X$.  

Let $\E$ be a complex vector bundle on $M$ equipped with a connection
$\nabla$. The pullback $\sigma^*\E$ of $\E$ to $X$ is then equipped with a
connection $\sigma^*\nabla$.

\begin{lemma} \label{_autodua_(1,1)-on-twi_Lemma_}
The connection $\nabla$ is autodual if and only if the connection
$\sigma^*\nabla$ has curvature of Hodge type $(1,1)$. 
\end{lemma}

\proof Indeed, the curvature $R_X$ of $\sigma^*\nabla$ is equal to the pullback
$\sigma^*R_M$ of the curvature $R_M$ of $\nabla$. Therefore it is of Hodge
type $(1,1)$ on $X$ if and only if for every $I \in \cp$ the form $R_M$ is
of type $(1,1)$ in the induced complex structure $I$. By
Lemma~\ref{primitive} this happens if and only if $R_M$ is
$\su$-invariant. 
\endproof

\hspace{-8pt} In particular, for every autodual bundle $\langle \E, \nabla
\rangle$ the $(0,1)$-com\-po\-nent ${\sigma^*\nabla}^{0,1}$ of the
connection $\sigma^*\nabla$ satisfies ${\sigma^*\nabla}^{0,1} \circ
{\sigma^*\nabla}^{0,1} = 0$ and defines a holomorphic structure on the
bundle $\sigma^*\E$.

\begin{defn}
The holomorphic bundle $\langle \sigma^*\E, {\sigma^*\nabla}^{0,1} \rangle$
is called {\emp the twistor transform} of the autodual bundle $\langle \E,
\nabla \rangle$. 
\end{defn}

\subsection{$\C P^1$-holomorphic bundles over twistor spaces}

The twistor transform is in fact invertible. To construct an inverse
transform, we begin with some results on differential forms on the twistor
space $X$.

The product decomposition $X = M \times \cp$ induces the decomposition
$\A^1(X) = \sigma^*\A^1(M) \oplus \pi^*\A^1(\cp)$ of the bundle $\A^1(M)$ of
$1$-forms.  By the definition of the complex structure on $X$, the
projection onto the subbundle of $(0,1)$-forms commutes with the projection
onto the bundle $\pi^*\A^1(\cp)$. Therefore a Dolbeault
differential
$$
\bar\6_\cp:\A^0(X) \to \pi^*\A^{0,1}(\cp)
$$
is well-defined. 

\begin{defn}
A {\emp $\cp$-holomorphic bundle} on $X$ is a complex vector bundle $\E$ on
$X$ equipped with an operator $\bar\6_\cp:\E \to \E \otimes
\pi^*\A^{0,1}(\cp)$ satisfying 
$$
\bar\6_\cp(fa) = \bar\6_\cp(f)a + f\bar\6_\cp(a)
$$
for a function $f$ and a local section $a$ of $\E$. 
\end{defn}

\begin{rem}
For any point $m \in M$ the restriction $\wt{m}^*\E$ of a $\cp$-holomrphic
bundle $\E$ to the horizontal twistor line $\wt{m}:\cp \to X$ is
holomorphic in the usual sense.
\end{rem}

A $\bar\6_\cp$-closed smooth section $a \in \Gamma(X,\E)$ of a
$\cp$-holomorphic bundle $\E$ will be called {\emp $\cp$-holomorphic}.
Tensor products and $\Hom$-bundles of $\cp$-holomorphic bundles are again
$\cp$-holomorphic.  A differential operator $f:\E_0 \to \E_1$ will be
called {\emp $\cp$-holomorphic} if $\bar\6 f(a)=0$ for every local
$\cp$-holomorphic section $a$ of the bundle $\E$.

For every complex vector bundle $\E$ on $M$ the bundle $\sigma^*\E$ on $X$
is canonically $\cp$-holomorphic.  For a $\cp$-holomorphic bundle $\E$ let
$\sigma_*\E$ be the sheaf on $M$ of $\cp$-holomorphic sections of $\E$.
The functors $\sigma_*$ and $\sigma^*$ are adjoint.

\begin{defn}\label{const}
A $\cp$-holomorphic bundle $\E$ on $X$ is called {\emp $\cp$-constant} if it is
isomorphic to $\sigma^*\E_1$ for a complex vector bundle $\E$ on $M$. 
\end{defn}

Since $\sigma_*\sigma^*\E_1 \cong \E_1$ for every complex bundle $\E_1$ on
$M$ , a bundle $\E$ on $X$ is $\cp$-constant if and only if the canonical
map $\E \to \sigma^*\sigma_*\E$ is an isomorphism. The functor $\sigma^*$
is therefore an equivalence between the category of complex vector bundles
on $M$ and the category of $\cp$-constant $\cp$-holomorphic bundles on $X$.

\begin{lemma}\label{constant}
A $\cp$-holomorphic bundle $\E$ on $X$ is $\cp$-constant if and only if for
every horizontal twistor line $\wt{m}:\cp \to X$ the restriction
$\wt{m}^*\E$ is trivial.
\end{lemma}

\proof Clear. \endproof 

Note that if $\dim\Gamma(\cp,\wt{m}^*\E)$ is the same for every horizontal
twistor line $\wt{m}:\cp \to X$, then $\sigma_*\E$ is the sheaf of smooth
sections of a complex vector bundle on $M$.

\subsection{Differential forms and $\C P^1$-holomorphic bundles}

Let $\A^\cdot_M(X) = \sigma^*\A^\cdot(M) \subset \A^\cdot(X)$ be the
subcomplex of relative $\C$-valued forms on $X$ over $\cp$ and let
$\A^{0,\cdot}_M(X)$ be the quotient complex of forms of Hodge type
$(0,\cdot)$. Denote by $d_M$ and $\bar\6_M$ the corresponding differentials
and let $P:\A^\cdot_M \to \A^{0,\cdot}_M$ be the natural projection. By
definition all the bundles $\A^\cdot_M$ and $\A^{0,\cdot}_M$ are
$\cp$-holomorphic, and the differentials $d_M$ and $\bar\6_M$ are
$\cp$-holomorphic.

Let $\E$ be a holomorphic bundle on $X$. Then the complex structure
operator $\bar\6:\E \to \E \otimes \A^{0,1}(X,\E)$ can be decomposed 
$$
\bar\6 = \bar\6_M + \bar\6_\cp
$$
into an operator $\6_\cp:\E \to \otimes \E \otimes \pi^*\A^{0,1}(\cp)$ and
an operator $\bar\6_M:\E \to \E \otimes \A^{0,1}_M$. The operator
$\bar\6_\cp$ is a $\cp$-holomorphic structure on the bundle $\E$, and the
operator $\bar\6_M$ is $\cp$-holomorphic. 

This constrution is in fact invertible. Namely, we have the following. 

\begin{lemma}\label{rel.holo}
The correspondence $\E \mapsto \langle \E, \bar\6_M\rangle$ is an
equivalence between the category of holomorphic bundles on $X$ and the
category of $\cp$-ho\-lo\-mor\-phic bundles $\E$ on $X$ equipped with a
$\cp$-holomorphic operator $\bar\6_M:\E \to \E \otimes \A^{0,1}_M$
satisfying $0 = \bar\6_M^2:\E \to \E \otimes \A^{0,2}_M$ and 
$$
\bar\6_M(fa) = \bar\6(f)a + f \bar\6_M(a)
$$
for a function $f$ and a local section $a$ of $\E$.
\end{lemma}

\proof Clear. \endproof 

In particlular, every holomorphic bundle $\E$ on $X$ is canonically
$\cp$-holomorphic. We will call a holomorphic bundle $\E$ on $X$ {\emp
$\cp$-constant} if the corresponding $\cp$-holomorphic bundle is
$\cp$-constant in the sense of Definition~\ref{const}.

\subsection{The complex 
${\cal A}^{\bullet}_{{\mathrm top}}(M)$: definition.}

For any horizontal twistor line $\wt{m}:\cp \to X$ the restriction
$\wt{m}^*\A^i_M$ is trivial, while $\wt{m}^*\A^{0,i}_M$ is a sum of several
copies of the bundle $\calo(i)$ on $\cp$ (see, e.g.,
\cite{HKLR}). Therefore $\sigma_*\A^*\cdot_M \cong \A^\cdot(M,\C)$, and the
map $P$ induces a projection
$$
P:\A^\cdot(M,\C) \cong \sigma_*\A^\cdot_M(X) \to \sigma_*\A^{0,1}_M(X).
$$
Let $\A^{1,0}_M(X)$ be the sheaf of relative forms of type $(1,0)$. By
definition we have an exact sequence 
$$
0 \to \A^{1,0}_M(X) \to \A^1_M(X) \to \A^{0,1}_M(X) \to 0
$$
of $\cp$-holomorphic bundles on $X$. Consider the associated long exact
sequence for $\sigma_*$. The restriction of $\A^{1,0}_M$ to any horizontal
twistor line $\wt{m}:\cp \to X$ is a sum of several copies of $\calo(-1)$.
Therefore this long sequence reduces to the map
$$
\sigma_*\A^1_M(X) \overset{P}{\longrightarrow} \sigma_*\A^{0,1}_M(X),
$$
which is therefore an isomorphism. 

Recall that the bundle $\A^i(M)$ carries a representation of the group
$\su$ for every $i \geq 0$. This representation is completely reducible
\footnote{$\su$ acts along the fibers, which are finite-dimensional},
and contains isotypical components of highest weights $\leq i$. Let
$\A^i_\topp \subset \A^i$ be the component of highest weight exactly $i$. 

\begin{lemma}\label{iso.complex}
The map $P:\A^\cdot(M,\C) \to \sigma_*\A^{0,\cdot}_M(X)$ is compatible with the
$\su$-action on $\A^\cdot(M,\C)$.  The restriction $P\A^\cdot_\topp(M,\C)
\to \sigma_*\A^{0,\cdot}_M$ is an isomorphism.  
\end{lemma}

\proof 
We already know the claim for $\sigma_*\A^{0,1}_M$. Let $i > 1$. Note that
$\sigma_*\A^{0,i}_M$ is equipped with an $\su$-action by the Borel-Weyl
theory. The corresponding representation is of highest weight $i$. The map
$P$ is obviously compatible with this action, which proves the first
statement. To prove the second, it is enough to prove that $P$ is
invertible on the subbundles of highest vectors. But both these subsundles
equal $\A^{0,i}(M)$.
\endproof 

\subsection{The complex ${{\cal A}}^{\bullet}_{{\mathrm top}}(M)$ and autodual bundles.}

The complex $\A^\cdot_\topp(M) \cong \sigma_*\A^{0,1}_M$ plays the same
role for autodual bundles as the Dolbeault resp. de Rham complexes play for
holomorphic resp. flat ones.  Precisely, let $\E$ be a complex vector
bundle on $M$ equipped with a connection $\nabla:\E \to \A^1(\E)$. Extend
the operator $\nabla$ to a differential operator $D:\A^\cdot_\topp(\E) \to
\A^{\cdot+1}_\topp(\E)$ by means of the embedding $\A^\cdot(\E)
\hookrightarrow \A^\cdot(\E)$ and the natural $\su$-invariant projection
$\A^{\cdot+1}(\E) \to \A^{\cdot+1}_\topp(\E)$.

\begin{lemma}\label{autodual.complex}
The connection $\nabla$ is autodual if and only if its extension $D$
satisfies $D^2 = 0$. 
\end{lemma}

\proof The operator $D^2$ is the multiplication by the $\A^2_\topp$-part of
the curvature $R$ of the bundle $\E$ with respect to the decomposition
$$
\A^2(\End\E) = \A^2_\topp(\End\E) \oplus \A^2_\inv(\End\E). 
$$
Thus it vanishes if and only if $R$ is $\su$-invariant, which by definition
means that $\nabla$ is autodual.
\endproof  

Let $\E$ be a complex vector bundle equipped with an autodual connection
$\nabla$. Since $\A^1(\E) \cong \sigma_*\A^{0,1}_M(X) \times \E$, the map
$\nabla:\E \to \E \otimes \sigma_*\A^{0,1}_M(X)$ defines a $\cp$-holomorphic
map 
$$
\bar\6_M:\sigma^*\E \to \sigma^*\E \otimes \A^{0,1}_M(X)
$$
of $\cp$-holomorphic bundles on $X$. By Lemmas~\ref{iso.complex} and
\ref{autodual.complex} the map $\bar\6_M$ extends to a map $\A^{0,\cdot}_M(X)
\otimes \sigma^*\E \to \A^{0,\cdot+1}_M(X) \otimes \E$ satisfying $\bar\6_M^2
= 0$. By Lemma~\ref{rel.holo} this map defines a holomorphic structure on
the bundle $\sigma^*\E$. 

\begin{lemma}
The holomorphic bundle $\langle\sigma^*\E,\bar\6_M\rangle$ on $X$ is 
isomorphic to the twistor transform of the autodual bundle $\E$. 
\end{lemma}

\proof 
Clear. 
\endproof 

Let now $\E$ be an arbitrary $\cp$-constant holomorphic bundle on $X$. 
Then the sheaf
$\sigma_*\E$ is the sheaf of sections of a vector bundle. Since $\E$ is
$\cp$-constant, $\E \cong \sigma^*\sigma_*\E$. 
The operator 
$$
\6_M:\E \cong \sigma^*\sigma_*\E \to \E \otimes \A^{0,1}_M(X)
$$
gives by adjuction an operator 
$$
\nabla:\sigma_*\E \to \sigma_*\left(\E \otimes \A^{0,1}_M(X)\right) \cong 
\sigma_*\E \otimes \A^1_\topp(M).
$$
By Lemmas~\ref{iso.complex} and \ref{autodual.complex} the operator
$\nabla$ is an autodual connection on $\sigma_*\E$.

\begin{defn}
The autodual bundle $\langle \sigma_*,\nabla \rangle$ on $M$ is called
{\emp the inverse twistor transform} of the $\cp$-constant holomorphic bundle
$\E$ on $X$.
\end{defn}

\begin{theorem} \label{_twisto_transfo_equiva_Theorem_}
The direct and inverse twistor transforms are mutually inverse equivalences
between the category of autodual bundles on $M$ and the category of
$\cp$-constant holomorphic bundles on $X$.
\end{theorem}

\proof 
Clear.
\endproof 

\section{Stability of the twistor transform.}
\label{_stabi_of_twi_tra_Section_}

\subsection{Introduction}
Let $M$ be a hyperk\"ahler manifold and let $X$ be its twistor
space. Consider a semistable autodual bundle $\langle \E,\nabla \rangle$ on
$M$ and let $\sigma^*\E$ be its twistor transform. The bundle $\sigma^*\E$
is a holomorphic bundle on $X$. In this section we prove under certain
conditions that $\sigma^*\E$ is semistable in the sense of
\ref{ss.twistor}. More precisely, we have the following.

\begin{prop}
Let $M$ be a hyperk\"ahler manifold. Denote its twistor space by $X$. Let
$\langle \E, \nabla \rangle$ be a semistable autodual bundle on $M$ and let
$\sigma^*\E$ be its twistor transform. Then for every coherent subsheaf $\F
\subset \sigma^*\E$ we have
$$
\frac{\deg c_1(\F)}{\rank\F} \leq \frac{\deg c_1(\sigma^*\E)}{\rank\sigma^*\E},
$$
where $\rank\F$ is the rank of the generic fiber of $\F$ and $\deg$ is
understood in the sense of subsection~\ref{ss.twistor}. 
\end{prop}

\subsection{Semistability for $\C P^1$-constant bundles}
Before we give a proof of this Proposition, we prove the following. 

\begin{lemma}\label{constant.ss}
Let $\E$ be a holomorphic bundle on the twistor space $X$. Assume that $\E$
is $\cp$-constant (that is, for every horizontal twistor line $\wt{m}:\cp
\to X$ the restriction $\wt{m}^*\E$ is trivial). Then $\deg \E = 0$, and the
bundle $\E$ is semistable.
\end{lemma}

\proof Since $H^1(\cp,\R) = 0$, 
\begin{equation}\label{ff}
H^2(X,\R) = H^2(M,\R) \oplus H^2(\cp,\R).
\end{equation}
Since $\E$ is $\cp$-constant, $c_1(\E) \in H^2(M,\R)$. For every $I \in
\cp$ let $X_I = \pi^{-1}(I) \subset X$ be the fiber over $I$. Since   
$$
c_1(\E)_{X_I} = c_1(\E|_{X_I}) \in H^{1,1}_I(M) 
$$
is of Hodge type $(1,1)$ for every $I \in \cp$, the class $c_1(\E) \in
H^2(M,\R)$ is $\su$-invariant by Lemma~\ref{primitive}. Therefore $\deg
c_1(\E) = \Lambda c_1(\E) = 0$.

To prove semistability, let $\F \subset \E$ be a coherent subsheaf. It is
enough to prove that $\deg c_1(\F) \leq 0$. Let $c_1(\F) = c_M + c_\cp$ be the
decomposition associated with \eqref{ff}. Again, by Lemma~\ref{primitive}
$c_M$ is $\su$-invariant, and $\deg c_1(\F) = \deg c_\cp$. For a generic
horizontal twistor line $\wt{m}:\cp \to X$ we have $c_\cp = c_1(\wt{m}^*\F)
\in H^2(\cp,\R)$. Since $\wt{m}^*\E$ is trivial, it is semistable, and
$\deg c_1(\wt{m}^*\F) \leq 0$.  \endproof

Let ${\cal M}^{ss}_X$ be the moduli space of semistable holomorphic bundles
on $X$. Lemma~\ref{constant.ss} implies that the set $\M_{const}$ of
equivalence classes of $\cp$-constant holomorphic bundles on $X$ is a
subset of ${\cal M}^{ss}_X$. The subset $\M_{const} \subset {\cal M}^{ss}_X$ is open.

\subsection{Conclusion}
The Propostion now follows directly from Lemma~\ref{constant.ss}. Moreover,
the twistor transform provides an isomorphism $\M_\inv \to \M_{const}
\subset {\cal M}^{ss}_X$ from the moduli space $\M_\inv$ of autodual
bundles on $M$ to the open subset of $\cp$-constant bundles in the moduli
space ${\cal M}^{ss}_X$

\section[Stable bundles and projective lines in twistor spaces.]
{Stable bundles and projective lines \\ in twistor spaces.}
\label{_lines_Section_}

\subsection{Hyperk\"ahler structure on the Mukai dual space}
\label{_Mukai_dual_Subsection_}

Let $M$ be a compact hyperk\"ahler manifold and let $\E$ be a complex
vector bundle on $M$ with $\su$-invariant Chern classes $c_1(\E)$ and
$c_2(\E)$. Consider the moduli space $\M_0$ of stable holomorphic
structures on $\E$ and let ${\cal M}^{reg}_0 \subset \M_0$ be the dense
open subset of smooth points in $\M_0$. 
Recall that the subset ${\cal
M}^{reg}_0$ is equipped with a natural K\"ahler metric, called {\em the
Weil-Peterson metric}.

It was proved in \cite{Vb} that the Weil-Peterson metric on ${\cal
M}^{reg}_0$ is actually hyperk\"ahler. Moreover, the complex manifold
$\left({\cal M}^{reg}_0\right)_J$ with the complex structure induced by a
quaternion $J \in \cp \subset \h$ was naturally identified with the subset
of smooth points in the moduli space of stable holomorphic structures on
$\E$ with respect to the complex structure $J$ on $M$.

Let $\X_{reg}$ be the twistor space of the hyperk\"ahler manifold ${\cal
M}^{reg}_0$. Consider the topological space $\X = \M_0 \times \cp$. We
have a natural embedding $\X_{reg} \subset X$. In \cite{Vb} the complex
structure on $\X_{reg}$ and the real structure $\iota:\X_{reg} \to
\X_{reg}$ were naturally extended to the whole of $\X$. The
complex-analytic space $\X$ is in general singular. However, the
fundamental Theorem~\ref{inv} still holds for the natural projection
$\pi:\X \to \cp$. We will call holomorphic sections $\cp \to \X$ of the
projection $\pi:\X \to \cp$ {\emp twistor lines} in $\X$. The space $\M_0$
is then naturally isomorphic to the subset of real twistor lines in
$\X$. 

These data define {\bf singular hy\-per\-k\"ah\-ler structure}
on  $\M_0$ (see \cite{Vb} for details). The space $\M_0$
with this hyperk\"ahler structure
is called {\bf Mukai dual} to $M$ (results of \cite{Vb}
generalise Mukai's work about duality of K3 surfaces). 
We must caution the reader that this version of Mukai duality
is not involutive, as the term ``dual'' might erroneously imply.

\subsection{Fiberwise stable bundles}

Let $X$ be the twistor space of the hyperk\"ahler manifold $M$. Let
$\M_{const} \subset {\cal M}^{ss}_X$ be the open subset of $\cp$-constant
holomorphic structures in the moduli space ${\cal M}^{ss}_X$ of semistable
holomorphic structures on the bundle $\sigma^*\E$. In the last section we
have identified $\M_{const}$ with the space $\M_\inv$ of $(0,1)$-stable autodual
connections on the bundle $\E$. 

In this section we will need still another notion of stability for
holomorphic bundles over $X$. 

\begin{defn}\label{fib.st}
Call a stable holomorphic structure $\bar\6$ on $\sigma^*\E$ {\emp fiberwise
stable} if for any $L \in \cp$ the restriction of $\langle \sigma^*\E,
\bar\6 \rangle$ to the fiber $X_L = \pi^{-1}(L) \subset X$ is stable. 
\end{defn}

Let $\M_\fib \subset {\cal M}^{ss}_X$ be the subset of fiberwise stable
holomorphic structures. The intersection $\left(\M_{const} \cap
\M_\fib\right) \subset {\cal M}^{ss}_X$ is, then, isomorphic to the moduli
space of autodual connections on $\E$ inducing a stable holomorphic
structure on $\E$ for every $I \in \cp$.

The goal of this section is to prove the following.

\begin{theorem}\label{iso}
The space $\M_\fib$ is naturally isomorphic to the space $\Sec$ of twistor
lines in the manifold $\X$. 
\end{theorem}

\subsection{Stability of fiberwise-stable bundles}

We begin by noting that one of the conditions in Definition~\ref{fib.st} is
in fact redundant. 

\begin{lemma}\label{gen.st}
Let $\E$ be a holomorphic bundle on the twistor space $X$. If the
restriction $i^*\E$ is stable for a generic\footnote{In the sense of
\cite{Vsym}} point $I \in \cp$, then the bundle $\E$ is stable. 
\end{lemma}

\proof 
Indeed, it was proved in \cite{Vsym} that for a generic point $I \in \cp$
every rational $(1,1)$-cohomology class for the fiber $X_I$ is of degree
zero. Therefore a stable holomorphic bundle on $X_I$ has no proper
subsheaves. Hence for a proper subsheaf $\F \subset \E$ either $\F$ or
$\E/\F$ is supported on non-generic fibers of $\pi:X \to \cp$. In
particular, either $\F$ or $\E/\F$ is a torsion sheaf. This implies that
the bundle $\E$ is stable.
\endproof 

\subsection{Modular interpretation of the 
Mukai dual twistor space}

We now construct a map $\M_\fib \to \Sec$. To do this, we give a
modular interpretation of the space $\X$. 

For any point $I \in \cp$ let $i:X_I \hookrightarrow X$ be the natural
embedding of the fiber $X_I = \pi^{-1}(I) \subset X$.  For a stable
holomorphic bundle $\E$ on the fiber $X_I$ call the coherent sheaf $i_*\E$
on $X$ {\emp a stable sheaf on $X$ supported in $I$}, or simply a {\emp
fiber-supported stable sheaf}. 

More generally, for a complex analytic space $Z$ call a coherent sheaf
$\EE$ on $Z \times X$ {\emp a family of fiber-supported stable sheaves on
$X$} if there exists a holomorphic map $f_Z:Z \to \cp$ such that $\E$ and 
a holomorphic bundle $\EE_0$ on the subspace $Z \times_\cp X \subset Z
\times X$ such that 
\begin{enumerate}
\item $\EE \cong i_*\EE_0$, where $i:Z \times_\cp X \to Z \times X$ is the
natural embedding. 
\item For every point $z \in Z$ the restriction of $\EE$ to $z \times
\pi^{-1}(f_Z(z)) \subset Z \times X$ is a stable holomorphic bundle.
\end{enumerate}
The space $\X$ is obviously the moduli space for families of
fiber-supported stable sheaves on $X$. The holomorphic map $f_\X:\X \to
\cp$ is the natural projection.

Let now $\E$ be a fiberwise stable holomorphic bundle on $X$. For every $I
\in \cp$ the coherent sheaf $i_*i^*\E$ on $X$ is a stable sheaf supported
in $I$. The correspondence $I \mapsto i_*i^*\E$ defines a holomorphic map
$\cp \to \X$. This map is a section of the projection $\X \to \cp$, hence
defines a point $\psi(\E) \in \Sec$. The correspondence $\E \mapsto
\psi(\E)$ comes from a holomorphic map $\psi:\M_\fib \to \Sec$ of the
corresponding moduli spaces. 

\subsection{Coarse and fine moduli spaces: a digression}

In order to prove that the map $\psi:\M_\fib \to \Sec$ is an isomorphism,
we need to make a digression about universal objects and coarse moduli
spaces. 

Let $\Var$ be the category of complex-analytic varieties and let $\F:\Var
\to \Sets$ be a functor. Recall that a complex-analytic space $Y$ is said to
be a fine moduli space for the functor $\F$ if $\F \cong \Hom(\bullet, Y)$.
This implies that there exists an element $C \in \F(Y)$ such that for every
complex-analytic space $U$ and an element $a \in \F(U)$ there exists a
unique map $f:U \to Y$ such that $a =\F(f)(C)$. Such an element $C$ is
called {\em the universal solution} to the moduli problem posed by $\F$.

It is well-known that geometric moduli problems only rarely admit fine
moduli spaces. The common way to deal with this is to introduce a weaker
notion of a {\em coarse moduli space}. For the purposes of this paper the
following notion suffices. 

\begin{defn}
A complex-analytic space $Y$ is called a {\emp coarse moduli
space for the problem posed by $\F$} if for any complex-analytic space $Z$
and an element $a \in \F(Z)$ there exist a unique map $f:Z \to Y$, an
open covering $U_\alpha$ of the space $Y$ and a collection $C_\alpha \in
\F(U_\alpha)$ such that for every index $\alpha$
$$
\F(f)(C_\alpha) = a|_{f^{-1}(U_\alpha)}.
$$
\end{defn}

Heuristically, a coarse moduli space $Y$ admits locally a universal solution
for the moduli problem $\F$, but these solutions need not come from a single
global solution in $\F(Y)$. 

All the moduli spaces constructed as infinite-dimensional quotients by
means of the slice theorem are coarse moduli spaces in the sense of this
definition. This applies to all the moduli spaces considered in this paper,
and to the space $\X$ in particular. Therefore for every point $x \in \X$
there exists a neighborhood $U_x \subset X$ and a coherent sheaf $\EE_x$ on
$U \times X$ which is a family of fiber-supported stable sheaves on $X$
universal for the moduli problem. Let $U$ be such a neighborhood. Then the
universality of the sheaf $\EE$ implies that
$$
\Aut\EE = \Gamma(U,\calo^*). 
$$

\begin{lemma}\label{coarse}
Let $\cp \to \X$ be a section of the projection $\pi:\X \to \cp$. There
exists a coherent sheaf\/ $\EE$ on $\cp \times X$ such that for every $x
\in \cp \subset \X$ the restriction\/ $\EE|_{U_x \times X}$ is isomorphic
to the universal sheaf\/ $\EE_x$.
\end{lemma}

\proof Indeed, cover $\cp$ by open subsets of the form $U_x$ and choose a
finite subcovering $U_\alpha$. In order to define a sheaf $\EE$, it is
enough to choose a system of isomorphisms 
$$
g_{\alpha\beta}:\EE_\alpha|_{(U_\alpha \cap U_\beta) \times X} \to
\EE_\beta|_{(U_\alpha \cap U_\beta) \times X} 
$$
for every intersection $U_\alpha \cap U_\beta$ so that $g_{\alpha\beta}
\circ g_{\beta\gamma} = g_{\alpha\gamma}$ for every three indices
$\alpha,\beta,\gamma$. Since $\Aut\EE_\alpha = \calo^*_{U_\alpha}$, the 
obstruction to finding such a system of isomorphisms lies in the second
\v{C}ech cohomology group $H^2(\cp,\calo^*)$. Consider the long exact
sequence
$$
H^2(\cp ,\calo) \longrightarrow H^2(\cp, \calo^*) \longrightarrow H^3(\cp,
\Z)
$$
associated to the exponential exact sequence
$$ 
0 \longrightarrow \Z \longrightarrow \calo \longrightarrow \calo^*
\longrightarrow 0.
$$ 
Since $H^2(\cp, \calo) = H^3(\cp, \Z) = 0$, the group $H^2(\cp ,\calo^*)$
vanishes.  
\endproof

\subsection{Conclusion}
We can now finish the proof of Theorem~\ref{iso}. It remains to prove that
the map $\psi:\M_\fib \to \X$ is an isomorphism.  We will construct an
inverse map $\psi^{-1}:\X \to \M_\fib$.

Let $x \in \Sec$ be a point and let $\wt{x}:\cp \to \X$ be the
corresponding section.  Let $\EE$ be the coherent sheaf on $\wt{x}(\cp)
\times X$ constructed in Lemma~\ref{coarse}.  Let $\Delta = \pi \times
\id:X \hookrightarrow \cp \times X$ be the embedding of $X$ into $\cp
\times X$ as the preimage of the diagonal under the natural projection $\id
\times \pi:\cp \times \X \to \cp \times \cp$. The sheaf $\EE$ is by
definition isomorphic to the direct image of a
holomorphic vector bundle $\E$ on $X$: $\EE \cong \Delta_*\E$.

For every point $I \in \cp$ the coherent sheaf $i_*i^*\E$ on $X$ is
canonically isomorphic to the restriction of $\EE$ to $I \times X \subset
\cp \times X$. Therefore the bundle $\E$ is stable by
Lemma~\ref{gen.st}. Let $\psi^{-1}(x) \in \M_\fib$ be the corresponding
point in the moduli space $\M_\fib$.

By construction $\psi(\psi^{-1}(x)) = x$. To prove that $\psi^{-1} \circ
\psi = \id$, consider a stable bundle $\E \in \M_\fib$. Let $\EE$ be the
coherent sheaf on $\cp \times X$ constructed in Lemma~\ref{coarse} and let
$p:\cp \times X \to X$ be the projection onto the second factor. By
definition
$$
\EE \cong \Delta_*\Delta^*(p^*\E) \cong \Delta_*(\Delta \circ p)^*\E \cong
\Delta_*\E. 
$$
Therefore $\psi^{-1}(\psi(\E)) = \E \in \M_\fib$. This finishes the proof
of Theorem~\ref{iso}.

\section{Conjectures and open questions.}
\label{_conje_Section_}

\subsection{NHYM moduli spaces and hyperk\"ahler reduction}
\label{_Hyperkae_redu_Subsection_}

\subsubsection{}
Let $M$ be a K\"ahler manifold and let $\M$ be the moduli space of NHYM
connections on a complex bundle $\E$ over $M$. We have shown in
Section~\ref{_NHYM_Section_} 
%
%
that the space $\M$ is equipped with a natural closed holomorphic $2$-form
$\Omega$ which is is symplectic at least in a neighborhood of the subset of
Hermitian connections. In fact one could hope for a much stronger
statement. 

\begin{conjecture}\label{hyp.nhymspace}
There exists a hyperk\"ahler metric on $\M$ such that $\Omega$ is the
associated holomorphic symplectic from. 
\end{conjecture} 

Note that the construction of the form $\Omega$ is completely parallel to a
construction of a holomorphic symplectic form on the Hitchin-Simpson moduli
space $\M_{DR}$ of flat connections on $\E$ (\cite{S2}). The analog of
Conjecture~\ref{hyp.nhymspace} for $\M_{DR}$ is known. 

\subsubsection{}
To provide some evidence for Conjecture~\ref{hyp.nhymspace}, we give an
interpretation of the NHYM equation in the context of hyperk\"ahler
reduction. 

Let $\A$ be the space of all connections on the complex vector
bundle $\E$. The space $\A$ is an affine space over the complex vector
space $\A^1(M,\End\E)$ of $\End\E$-valued $1$-forms on $M$. Choose an
Hermitian metric $h$ on the bundle $\E$. The decomposition 
$$
\A^1(M,\End\E) = \A^{1,0}(M,\End\E) \oplus \A^{0,1}(M,\End\E)
$$
allows one to define a quaternionic structure on the space
$\A^1(M,\End\E)$. Together with the natural trace metric, this structure
makes the space $\A$ an (infinite-dimesional) hyperk\"ahler manifold. 

The complex gauge group $\G = \Maps(M,\Aut\E)$ acts on the space $\A$. This
action is compatible with the hyperk\"ahler structure on $\A$. Therefore one
can apply to the space $\A$ the machinery of {\emp hyperk\"ahler reduction}
(see \cite{HKLR}). It turns out that 
the complex moment map  $\A \to \C$ is equal to the map 
$YM:\A \to \Gamma(M,\End\E)$, $\nabla \mapsto
\Lambda\nabla^2$. Vanishing of this map is precisely
the NHYM condition.

Let $\A_0 = YM^{-1}(0) \subset \A$ be the subset of connections with
$\Lambda\nabla^2 = 0$. By the general principles of hyperk\"ahler reduction
the quotient $\A_0 / \G$ should be hyperk\"ahler. The NHYM moduli space
$\M$ is the closed subset $\M \subset \A_0 / \G$ of equivalence classes of
connections with curvature $R=\nabla^2$ which satisfy $\Lambda R = 0$ and
are, in addition, of Hodge type $(1,1)$. We expect that the embedding $\M
\hookrightarrow \A_0 / \G$ is compatible with the hyperk\"ahler structure
on $\A_0 /\ G$ and gives a hyperk\"ahler structure on $\M$ by restriction.

\subsubsection{}
The hyperk\"ahler reduction construction of the NHYM moduli space also
allows to formulate an analog of the Uhlenbeck-Yau Theorem for
NHYM-bundles. We first give a new definition of stability for NHYM bundles,
more natural than the $(0,1)$-stability used in the body of the paper.

Let $\overline{M}$ be the complex-conjugate complex manifold to $M$. Since
$M$ and $\overline{M}$ are the same as smooth manifolds, the bundle $\E$
can be also considered as a complex vector bundle on $\overline{M}$. For every
connection $\nabla$ on $\E$ the $(1,0)$-part $\nabla^{1,0}$ defines a
holomorphic structure on the complex bundle $\E$ on $\overline{M}$. 

\begin{defn}
Let $\langle\E, \nabla\rangle$ be a bundle with a $(1,1)$-connection over a
complex manifold $X$. Let $U\subset X$ be a Zariski open subset in $X$, and
let $\F\subset \E\restrict{U}$ be a subbundle which is preserved by
$\nabla$. Then $\F$ is called {\bf a subsheaf of $\langle\E, \nabla\rangle$} if
the following two conditions hold:

\begin{description}
\item[(i)] Consider $\E$ as a holomorphic bundle over $M$,
with a holomorphic structure defined by the $(0,1)$-part 
of the connection. Then there exist a coherent subsheaf
$\widetilde \F\subset \E$ on $M$ such that the restriction
$\widetilde \F\restrict{U}$ is a sub-bundle of $\E$ which 
coinsides with $\F$.

\item[(ii)] 
Consider $\E$ as a holomorphic bundle over $\overline{M}$,
with a holomorphic structure defined by the $(1,0)$-part 
of the connection. Then there exist a coherent subsheaf
$\widetilde \F\subset \E$ on $\overline{M}$ such that the restriction
$\widetilde \F\restrict{U}$ is a sub-bundle of $\E$ which 
coinsides with $\F$.

\end{description}
\end{defn}

For a subsheaf $\F\subset \E$, it is straightforward to define 
the Chern classes and the degree. As usually, $\F$ is called
{\bf destabilizing} if 

\[ \frac{\deg \F}{\rank \F} 
   \geq \frac{\deg {\E}}{\rank {\E}} 
\]

\begin{defn}\label{_stabili_for_hyper_redu-Definition_}
Let $\langle\E, \nabla\rangle$ be a bundle with $(1,1)$-connection over a 
compact K\"ahler manifold $M$. Then $\langle\E, \nabla\rangle$ is called
$\nabla$-stable if there are no destabilizing subsheaves 
$\F\subset \langle\E, \nabla\rangle$.  
\end{defn}

This definition generalizes the definition \cite{_Simpson:harmonic_}
of stability for flat bundles.

\begin{rem}
Clearly, for NHYM bundles, $(0,1)$-stability implies the stability in the
sense of Definition \ref{_stabili_for_hyper_redu-Definition_}
\end{rem}

An analogy with the Kempf-Ness Theorem suggests that every stable
$\G$-orbit in $\A_0$ has non-trivial intersection with the zero set of the
real moment map $\A_0 \to \Gamma(M,\End_\R\E)$ from $\A_0$ to the space of
anti-Hermitian endomorphisms of the bundle $\E$. This moment map can be
described more explicitly. 

\begin{defn}[pseudocurvature]
Let $\langle\E, \nabla\rangle$ be a bundle with $(1,1)$-connection and a
Hermitian metric $h$, not necessary compatible.  
Let $\nabla = \nabla' + \nabla''$ be the decomposition
of $\nabla$ onto $(1,0)$ and $(0,1)$-parts. Consider the connection in
$\overline{\E}^*$ associated with $\nabla$. Since $h$ identifies $\E$ and
$\overline{\E}^*$, this gives another connection in $\E$, denoted by
$\nabla_h$.  The average $\nabla_h= \frac{\nabla +\nabla_h}{2}$ is again a
connection, and is compatible with $h$. Let $\theta$ be the difference
$\theta:= \frac{\nabla -\nabla_h}{2}$, which is a tensor.  Applying
$\nabla_h$ to $\theta$, we obtain a 2-form $\Xi$ with coefficients in
$\End(\E)$. This form $\Xi$ is called {\bf the pseudocurvature} of the
triple $\langle \E, \nabla, h\rangle$.
\end{defn}

It turns out that the real moment map on a NHYM connection $\nabla$ 
is given by 
$$
\nabla \mapsto \Lambda(\Xi),
$$ 
where $\Xi$ is the pseudocurvature. 

\begin{defn} \label{_harmonic_me_Definition_}
Let $\langle\E, \nabla\rangle$ be a bundle with a NHYM $(1,1)$-connection,
and let $h$ be an Hermitian metric on $\E$, not necessarily compatible with
$\nabla$ Then $h$ is called {\bf harmonic} if $\Lambda \Xi =0$, where $\Xi$
is the pseudocurvature of $\langle \E, \nabla, h \rangle$.
\end{defn}

\begin{conjecture} 
Let $\langle\E, \nabla\rangle$ be a bundle with NHYM connection $\nabla$. 
Then there exists a harmonic metric $h$ on $\E$ if and only if $\E$ is a
direct sum of $\nabla$-stable bundles.  Also, if $\E$ itself is
$\nabla$-stable, then $h$ is unique, up to a constant factor.
\end{conjecture}

An analogous statement is known for flat connections. See
\cite{_Simpson:harmonic_} for a discussion. 

\subsection{K\"ahler base manifold: open questions.}

In this subsection, we relate questions pertaining to the
case of base manifold $M$ compact and K\"ahler, but not necessarily
hyperk\"ahler.

\begin{question} \label{_NHYM_flat?_Question_}
Let $(B, \nabla)$ be a NHYM-connection in a bundle with zero 
Chern classes. Is it true that $\nabla$ is necessarily flat?
\end{question}

In Hermitian case, the answer is affirmative 
by L\"ubcke \cite{_Lubcke_}
and Simpson \cite{S}. In a neighbourhood of Hermitian Yang-Mills
connection, all NHYM connections on a bundle with zero 
Chern classes are also flat,
at one can see, e. g., from Proposition \ref{series}.
The hyperk\"ahler analogue of this question is Question 
\ref{_NHYM-are-autodu_Question_}.

\hfill

Let $B$ be a stable holomorphic bundle over $M$ and let $St(B)$ be the
deformation space of stable holomorphic structures on $B$.
In Section \ref{_NHYM_Section_}, we defined a Kuranishi map
$\phi:\; U \hookrightarrow H^1(\End(B))$, where $U$ is a neighbourhood
of $[B]$ in $St(B)$. The map $\phi$ is, locally, a closed embedding,
and its image in a neighbourhood of zero in $H^1(\End(B))$ is an algebraic
subvariety, defined by the zeroes of so-called Massey products.
Let $C$ be the Zariski closure of the image 
$\phi(U)$ in $H^1(\End(B))$. Let $\operatorname{NHYM}(B)$ be the space
of NHYM connections inducing the same holomorphic structure.
In \eqref{_from_NHYM_to_classes_Equation_}, we construct
the map $\operatorname{NHYM}(B) \stackrel{\rho}{\arrow} \overline C$,
where $\overline C$ is a complex conjugate manifold to $C$,
and prove that in a neighbourhood of zero $\rho$ is isomorphism.

Two questions arise:

\begin{question}\label{_rho_surjec_Question_}
Is the map $\rho$ surjective?
\end{question}

\begin{question} \label{_rho_etale_Question_}
Is the map $\rho$ etale? Bijective?
\end{question}

These two questions might be reformulated in a purely algebraic
way. Let $B$ be a stable holomorphic bundle equipped with a Hermitian
Yang-Mills metric, and $X$ be the space of all $(1,0)$-forms
$\theta\in \Lambda^{1,0}(\End (B))$ satisfying

\begin{equation} 
  \begin{cases}
	\6 \theta &= \theta\wedge \theta \\
	\6^* \theta& =0,
  \end{cases} 
\end{equation}
where $\6$ is the $(1,0)$-part of the connection, and 
$\6^*$ the adjoint operator. The middle cohomology space of the complex

\[ \Lambda^{2,0}(\End(B)) \stackrel{\6^*}{\arrow}
   \Lambda^{1,0}(\End(B)) \stackrel{\6^*}{\arrow}
   \End(B)
\]
is naturally isomorphic to the complex conjugate space  
to $H^1(\End(B))$. This gives a map
\[ 
    X \stackrel{\rho} {\arrow} \overline{H^1(\End(B))}, 
\]
associating to $\theta$ its cohomology class. Then,
Proposition \ref{kur} implies that the image of $\rho$ 
lies in $\overline C$ and locally in a neighbourhood of zero,
$\rho: \; X \arrow \overline C$ is an isomorphism.
Question \ref{_rho_surjec_Question_} asks 
whether $\rho$ is surjective onto
$\overline C$, and Question \ref{_rho_etale_Question_} 
asks whether $\overline \pi$ is etale,
or even invertible.

\subsection{Autodual and NHYM connections over a hyperk\"ahler base.}

\subsubsection{}
The first and foremost question (partially answered in 
Theorem \ref{_NHYM-are-autodu_Theorem_}; see also 
Question \ref{_NHYM_flat?_Question_}):

\begin{question} \label{_NHYM-are-autodu_Question_}
Let $(B, \nabla)$ be a NHYM bundle over a hyperk\"ahler
base manifold. Is $(B, \nabla)$ necessarily autodual?
\end{question}

\subsubsection{}
\label{compa_smoo_Mukai_du_Subsubsection_}
Let $M$ be a compact hyperk\"ahler manifold, and let $\cal S$ be
a connected component of the moduli of autodual connections
on a complex vector bundle ${\cal B}$. 
Assume that $\cal S$ contains a point $B$ 
which is Hermitian autodual. Consider the ``Mukai dual''
space $\widehat M$, that is, the moduli space 
of Hermitian autodual connections
on ${\cal B}$ (Subsection \ref{_Mukai_dual_Subsection_}). 
Assume that the connected component of $\widehat M$
containing $B$ is smooth and compact. Clearly, then, all connections
from $\cal S$ are fiberwise stable, in the sense of Definition
\ref{fib.st}. Thus, Theorem \ref{iso} gives an isomorphism
between $\cal S$ and the space $Sec(\widehat M)$ 
of twistor lines in $\Tw(\widehat M)$. 

In such situation, we are going to give a
conjectural description of the space $\cal S$, assuming
that the answer to 
\ref{_rho_surjec_Question_}---\ref{_rho_etale_Question_}
is affirmative.

\subsubsection{}
\begin{defn}[Twisted cotangent bundle] 
Let $M$ be a K\"ahler manifold, $\Omega^1M$ its
holomorphic cotangent bundle. The K\"ahler class 
\[ 
   \omega \in H^1(\Omega^1 M)= Ext^1(\calo (M), \Omega^1 M)
\] 
gives by Yoneda an exact sequence
\[
   0 \arrow \Omega ^1 M \arrow E \stackrel e \arrow \calo(M) \arrow 0, 
\]
where $\calo(M)$ is the trivial one-dimensional bundle. 
Let $\nu$ be a non-zero section of $\calo(M)$, and $E_\nu$
be the set of all vectors 
\[
   \left\{ v\in E\restrict m \;\;
   \left|\vphantom{\bigcup\limits_M^M}\right.
\;\; e(v) = \nu\restrict m\right\}
\]
where $m$ runs through all points of $M$. Consider $E_\nu$ as a
submanifold in the total space of $E$. Then $E_\nu$ is called
{\bf a twisted cotangent bundle of $M$}, denoted by $\Omega_\omega M$.
\end{defn}

The space $\Omega_\omega M$ has a natural action of $\Omega^1 M$
considered as a group scheme over $M$, and as such is a torsor
over $\Omega ^1 M$

\subsubsection{}
The affirmative answer to the stronger form
of \ref{_rho_etale_Question_} would give the proof
of the following conjecture.

\begin{conjecture} \label{_sec_to_twi_cota_Conjecture_}
Under assumptions of \ref{compa_smoo_Mukai_du_Subsubsection_},
there exists a natural isomorphism of complex manifolds
\[ Sec(\widehat M) \cong \Omega_\omega M, \]
where $\Omega_\omega M$ is the twisted cotangent bundle.
\end{conjecture}

In the general situation, there is a natural
map from the space of twistor lines $Sec(M)$ 
of a compact hyperk\"ahler manifold to $\Omega_\omega M$.
However, in general there are no approaches to the proof
of surjectivity.

\begin{example}
Let $M$ be a compact complex torus, $\dim_\C M = 2n$, and 
$B$ a trivial line bundle. Clearly, $M$ is hyperk\"ahler.
Then $\widehat M$ is  the dual torus, and $\cal S$ is the 
space of local systems on $M$, which is isomorphic to 
$(\C^*)^{2n}$. The space $C$ of \ref{_rho_etale_Question_}
is isomorphic to $H^1(\calo M)$, and the answer to
\ref{_rho_etale_Question_} is obviously affirmative. Thus,
$Sec(M)$ and $\Omega_\omega M$ are also isomorphic to $(\C^*)^{2n}$
and are Stein. 
\end{example}

In the following subsection, we shall see
that this is indeed a general phenomenon -- 
the space of twistor lines is equipped
with a canonical plurisubharmonic function and is
likely to be Stein. However, we don't know a general argument
constructing plurisubharmonic functions on the twisted
cotangent bundle -- this is one more mystery.

\subsection{Plurisubharmonic functions on moduli spaces.}

Let $M$ be a hyperk\"ahler manifold $\Tw \stackrel \pi \arrow \C P^1$
its twistor space, and $Sec(M)$ the space of sections
$s:\; \C P^1 \arrow \Tw$ of the map $\pi$, also called {\it twistor
lines} (Section \ref{_twistors_Section_}). It is easy to equip
$Sec(M)$ with a natural plurisubharmonic function.

Recall that $\Tw$ is isomorphic as a $C^\infty$-manifold to
$M \times \C P^1$. This decomposition gives a natural 
(non-K\"ahler) Hermitian metric on $\Tw$.

\begin{prop} \label{_volume_twi_line_plurisubharmo_Proposition_}
Consider the function $v:\; Sec(M) \arrow \R ^+$
which maps a line $s\in \Tw$ to its Hermitian volume, taken 
with respect to the Hermitian metric on $\Tw$. Then $v$
is strictly plurisubharmonic.
\end{prop}

\proof
Let $\omega$ be the differential 2-form which is 
the symplectic part of the Hermitian metric on $\Tw$.
Since the twistor lines are complex subvarieties in $\Tw$,
$v(s) = \int_s \omega$ for all twistor lines $s$. Then, for
all bivectors $x, \bar x$ in $T_s Sec(M)$, we have
\begin{equation}\label{_6bar6omega_Equation_}
  \6\bar\6 v(x, \bar x) = \int_s \6\bar\6 \omega({\bf x},\bar{\bf x}),
\end{equation}
where ${\bf x},\bar{\bf x}$ are the sections of $T \Tw\restrict{s}$
corresponding to $x, \bar x$.
Then, to prove that $v$ is plurisubharmonic it suffices
to show that $\6\bar\6 \omega({\bf x},\bar{\bf x})$ is positive. 
{}From Lemma \ref{_differe_of_Hermi_on_twistors_Lemma_}, it
is easy to see that 
$\6\bar\6 \omega = \omega\wedge \pi^* \operatorname{FS}(\C P^1)$,
where $\operatorname{FS}(\C P^1)$ is the Fubini-Study form on $\C P^1$. 
Clearly, then, $\6\bar\6\omega({\bf x},\bar{\bf x})$ is 
positive, and $v$ is plurisubharmonic. This proves Proposition 
\ref{_volume_twi_line_plurisubharmo_Proposition_}.
\endproof

One of the most intriguing questions of hyperk\"ahler geometry
is to learn whether the function $v$ is exhausting.

\hfill

In notation and assumptions of \ref{_sec_to_twi_cota_Conjecture_},
consider the space $Sec(\widehat M)$ which is isomorphic to
the space $\cal S$ of autodual connections. There is
the canonical Weil-Petersson metric on $\cal S$, coming
from results of Subsection \ref{_Hyperkae_redu_Subsection_}.
This metric is hyperk\"ahler. This metric is given by a potential,
which is equal to the integral of the square of the
absolute value of the curvature.

\begin{question}
Is the Weil--Petersson metric related to the metric given by $v$?
\end{question}

\subsubsection*{Acknowledgements:} The autors are grateful to
S.-T. Yau, who stimulated the interest to the problem, D. Kazhdan
and T. Pantev for valuable discussions, S. Arkhipov, M. Finkelberg 
and L. Positselsky for their attention, and to Soros Foundation 
which is our source of livelihood.

\end{document}